\documentclass[a4paper,fleqn]{cas-sc}

\pdfoutput=1
\ExplSyntaxOn
\keys_set:nn { stm / mktitle } { nologo }
\ExplSyntaxOff

\usepackage[square, numbers]{natbib}
\usepackage{algorithm,algpseudocode}

\usepackage{tikz}
\usepackage{amsmath}
\usepackage{rotating}
\usepackage{color,soul}
\usepackage{todonotes}
\usepackage[modulo]{lineno}
\usepackage{subfig}
\usepackage{caption}
\usepackage{arydshln}
\usepackage{schemata}
\newcommand\AB[2]{\schema{\schemabox{#1}}{\schemabox{#2}}}

\def\tsc#1{\csdef{#1}{\textsc{\lowercase{#1}}\xspace}}
\tsc{WGM}
\tsc{QE}
\tsc{EP}
\tsc{PMS}
\tsc{BEC}
\tsc{DE}

\begin{document}
\let\WriteBookmarks\relax
\def\floatpagepagefraction{1}
\def\textpagefraction{.001}
\shorttitle{A new mixed-integer programming model for irregular strip packing based on vertical slices with a reproducible survey}
\shortauthors{Lastra-Díaz\&Ortuño}

\title [mode = title]{A new mixed-integer programming model for irregular strip packing based on vertical slices with a reproducible survey}                      

\author[1]{Juan J. Lastra-Díaz}[orcid=0000-0003-2522-4222]
\cormark[1]
\cortext[cor1]{Corresponding author}
\ead{jlastra@ucm.es}

\author[1]{M. Teresa Ortuño}[orcid=0000-0002-5568-9496]
\ead{mteresa@ucm.es}

\address[1]{Department of Statistics and Operational Research, Institute of Interdisciplinary Mathematics, UCM Research Group HUMLOG, Complutense University of Madrid, Spain}

\begin{abstract}
The irregular strip-packing problem, also known as nesting or marker making, is defined as the automatic computation of a non-overlapping placement of a set of non-convex polygons onto a rectangular strip of fixed width and unbounded length, such that the strip length is minimized. Nesting methods based on heuristics are a mature technology, and currently, the only practical solution to this problem. However, recent performance gains of the Mixed-Integer Programming (MIP) solvers, together with the known limitations of the heuristics methods, have encouraged the exploration of exact optimization models for nesting during the last decade. Despite the research effort, the current family of exact MIP models for nesting cannot efficiently solve both large problem instances and instances containing polygons with complex geometries. In order to improve the efficiency of the current MIP models, this work introduces a new family of continuous MIP models based on a novel formulation of the NoFit-Polygon Covering Model (NFP-CM), called NFP-CM based on Vertical Slices (NFP-CM-VS). Our new family of MIP models is based on a new convex decomposition of the feasible space of relative placements between pieces into vertical slices, together with a new family of valid inequalities, symmetry breakings, and variable eliminations derived from the former convex decomposition. Our experiments show that our new NFP-CM-VS models outperform the current state-of-the-art MIP models. Finally, we provide a detailed reproducibility protocol and dataset based on our Java software library as supplementary material to allow the exact replication of our models, experiments, and results.
\end{abstract}




\maketitle

\section{Introduction}

Cutting and packing regular (convex) and irregular (non-convex) polygons onto a rectangular strip with unbounded length is a tedious and omnipresent task in most manufacturing industries based on the cutting of any flat material. For instance, \citet{Milenkovic1991-mo} study the nesting problem for the fashion and apparel industry, whilst \citet{Heistermann1995-qm} and \citet{Whelan1993-hf} do it on leather manufacturing, \citet{Elamvazuthi2009-pu} in furniture, \citet{Han2013-cb} in the glass industry, \citet{Alves2012-ab} in the automotive industry, and \citet{Cheok1991-jy} in shipbuilding. The irregular strip-packing methods aim to compute a non-overlapping placement of a set of irregular polygons onto a fixed-width rectangular strip with unbounded length, called the \emph{board}, whose length is the minimum between all feasible placements. Another closely related problem, called two-dimensional bin packing \cite{Lodi2002-jz}, is defined as the computation of a non-overlapping placement of a set of polygons onto a larger closed polygon, called the \emph{bin}, to minimize the number of bins required. Although most of bin packing problems are defined for rectangular bins and items \cite{Iori2021-bi}, we also find many irregular bin packing problems in the aforementioned industries. Strip and bin packing problems, and all their variants concerning the geometry of the pieces or boards, belong to the broader family of Cutting and Packing (C\&P) problems categorized by \citet{Dyckhoff1990-xv} and \citet{Wascher2007-fw}, and extensively reviewed by \citet{Sweeney1992-cb}, \citet{Dowsland1992-su}, \citet{Wang2002-sd}, and \citet{Bennell2013-jd}. 

Research on the irregular strip and bin packing problems can be traced back to the pioneering Linear Programming (LP) models for rectangular bin packing introduced by \citet{Gilmore1965-vd}, and the pioneering heuristic methods for irregular strip packing proposed by \citet{Art1966-vy}, \citet{Adamowicz1976-wv, Adamowicz1976-fp}, and \citet{Albano1980-tq} in the late nineteen sixties and seventies. These early works introduce many of the basic ideas subsequently exploited by all heuristics methods reported in the literature, such as the notion of a feasible non-overlapping region between pieces based on the No-Fit Polygon (NFP) representation, and the sequential placement of pieces based on a bottom-left heuristics. Given two polygons $A,B \subset \mathbb{R}^2$, the no-fit polygon of $B$ regarding $A$ is defined by $NFP_{AB} = A \oplus (-B(0,0))$, where $\oplus$ symbol denotes the Minkowski sum of two sets $S_1,S_2 \subset \mathbb{R}^2$, such that $S_1 \oplus S_2 = \{p + q: p \in S_1, q \in S_2 \}$ and $B(0,0)$ is the translation of the polygon $B$ such that its reference point is positioned at the origin. The boundary and outer region of the no-fit polygon $NFP_{AB}$ set the feasible positions in which polygon $B$ can be placed without rotating into a non-overlapping position regarding polygon $A$. The no-fit polygon allows computing the feasible relative placements for any polygon pair $(A, B)$ a priori by checking whether their relative position vector $(x_B - x_A, y_B-y_A)$ belongs to either the boundary or the outer region of $NFP_{AB}$. Consequently, this later property has converted the NFP into the most broadly adopted and effective geometric representation for the nesting problem reported in the literature, both by the families of heuristics methods \cite{Bennell2009-cy} and exact mathematical models \cite{Leao2020-bc}. 

\citet{Fowler1981-vr} and \citet{Milenkovic1991-mo} show that the irregular strip-packing problem is NP-complete. For this reason, most practical solutions reported in the literature since the pioneering work of \citet{Art1966-vy} are based on sequential placement heuristics to build efficiently feasible solutions that are combined with meta-heuristics for exploring the space of feasible solutions, as shown in most of surveys on nesting \cite{Dowsland1995-qp, Dowsland2002-hx, Hopper2001-rp, Bennell2009-cy, Riff2009-rv}. Current heuristics-based methods can efficiently compute acceptable solutions for large problem instances, which have encouraged their early and extensive adoption in all industries mentioned above, especially in those industries with a so significant and permanent diversity of products as the garment industry \cite{Puri2013-xo}. \citet{Elkeran2013-pe} introduces the current state-of-the-art heuristics method for nesting, called Guided Cuckoo Search (GCS), which defines a two-stages method based on piece clustering and a NFP-based bottom-left heuristics \cite{Gomes2002-qj} to build an initial feasible solution that is shrinked by solving an overlap minimization problem using a variant of the cuckoo search meta-heuristics. Although GCS was introduced almost a decade ago, subsequent works have been unable to outperform its results, as shown by \citet[table 3]{Pinheiro2016-bg}, \citet[table 5]{Sato2016-kc}, \citet[table 4]{Cherri2016-zr}, \citet[table 6]{Mundim2017-gs}, \citet[table 2]{Amaro_Junior2017-fi}, and \citet[tables 4-5]{Sato2019-lp}. GCS \cite{Elkeran2013-pe} also significantly outperforms a commercial nesting system broadly adopted in the industry \cite[figure 3]{Annamalai_Vasantha2016-na}. However, the heuristics methods demand several hours to improve their initial solutions \cite[table 6]{Mundim2017-gs}, and they can neither provide an optimality proof nor a gap measure to the optimal solution.

Recent performance gains of the Mixed-Integer Programming (MIP) solvers, and the limitations of the heuristics methods mentioned above, have encouraged the exploration of exact mathematical models for nesting during the last decade, as shown by \citet{Leao2020-bc}. The mathematical models for nesting can be categorized into continuous, discrete, or semi-continuous models according to the type of decision variable used to represent the position of the pieces. \citet[\S8]{Li1994-co} introduces the first continuous MIP model for nesting reported in the literature to solve a limitation of the pioneering LP compaction model of \citet{Li1993-uj, Li1995-aj}. Despite \citeauthor{Li1994-co}'s model is not experimentally evaluated, it sets the two main features of the family of continuous MIP models based on the NFP as follows: (1) the convex decomposition of the outer NFP feasible regions; and (2) the definition of mutually-exclusive binary variables to set the pairwise non-overlapping constraints between pieces defining the feasible regions for their relative placement. Thus, we say that \citeauthor{Li1994-co}'s model sets the basic continuous linear MIP model for nesting without rotations, from which all subsequent works, including the present work, introduce more tightened and refined formulations. For instance, \citet[\S5]{Dean2002-wx} introduces and evaluates for the first time a MIP model for nesting based on a refinement of the \citet{Daniels1994-zn} bin packing model that is essentially identical to the \citeauthor{Li1994-co}'s model. Subsequently, \citet{Fischetti2009-oa} (F\&L) introduce a refinement of \citeauthor{Li1994-co}'s model based on lifting the big-M formulation and a branching-priority algorithm to guide the Branch\&Bound (B\&B) exploration by fixing large sets of feasible relative placements among pieces. \citet{Alvarez-Valdes2009-yz} improve the F\&L model by introducing a new MIP model called HS2, which is based on a detailed convex decomposition of the feasible regions, a lifting for the bound constraints of the continuous variables, six new branching strategies, and the x-axis ordering of identical pieces to remove all symmetric solutions derived from their permutation. However, the HS2 model \cite{Alvarez-Valdes2009-yz} cannot significantly improve the F\&L model despite all their improvements.

More recently, \citet{Cherri2016-jf} introduce two continuous MIP models improving the HS2 model \cite{Alvarez-Valdes2009-yz}, together with the first MIP model integrating discrete rotations, which are based on the convex decomposition of the pieces and the definition of the non-overlapping constraints between pieces by using the convex no-fit polygons among their convex parts. The first model, called Direct Trigonometry Mode (DTM), uses a direct trigonometry function encoding the separation lines defined by the edges of the polygons to build the non-overlapping constraints between pairs of convex parts from two pieces, whilst the second model, called NoFit-Polygon Covering Model (NFP-CM), uses the convex NFP between convex parts and several valid inequalities. DTM and NFP-CM use the same x-axis symmetry-breaking for identical pieces proposed by \citet{Alvarez-Valdes2013-wg}. Subsequently, \citet{Rodrigues2017-dy} improve the NFP-CM model by breaking the symmetries of the feasible space for the relative placements between pieces, setting the current state-of-the-art in terms of performance among the family of exact continuous mathematical models for irregular strip packing. However, the Improved NFP-CM model \cite{Rodrigues2017-dy} is only able to solve small problem instances with up to 17 pieces with simple geometry \cite[table 2]{Rodrigues2017-dy}, and unlike the NFP-CM model \cite[table 3]{Cherri2016-jf}, it has not been evaluated on complex small and large problem instances yet.

The main aim of this work is to introduce a new family of continuous MIP models for irregular strip packing without rotations capable of improving the performance of current state-of-the-art MIP models. The new family of MIP models, called NFP-CM based on Vertical Slices (NFP-CM-VS), is a novel and more tightened formulation of the state-of-the-art NFP-CM models \cite{Cherri2016-jf, Rodrigues2017-dy} that is based on a new disjoint convex decomposition of the feasible regions between convex parts into vertical slices, together with a new family of valid inequalities, symmetry breakings, and variable reductions derived from the former geometric decomposition. A second aim of this work is to carry out a fair reproducible comparison of our new family of MIP models with the state-of-the-art family of NFP-CM models \cite{Cherri2016-jf, Rodrigues2017-dy} by replicating these later models and implementing our new MIP models into a single software platform integrating the latest version of the Gurobi solver, which will be provided together with a detailed reproducibility protocol and dataset as supplementary material to allow the exact replication of all our models, experiments, and results. 

\subsection{Motivation and hypothesis}
\label{sec:motivation}

Our main motivation is the proposal and evaluation of a new and more tightened formulation of the family of NFP-CM models introduced by \citet{Cherri2016-jf} and \citet{Rodrigues2017-dy}. Our main hypothesis is that the new NFP-CM-VS models could improve the performance of the current state-of-the-art MIP models.

The second motivation of this work is to implement a reproducible experimental survey for a fair and confirmatory comparison of all models, which is based on our software implementation of all MIP models evaluated herein into a single Java software library based on the latest version of the Gurobi solver, called Gurobi 9.5.

A third motivation is to bridge the lack of reproducibility resources hampering the independent replication and confirmation of previously reported models and results, as well as the incorporation of newcomers into this area, by providing a detailed reproducibility protocol and dataset as supplementary material to allow the exact replication of all MIP models evaluated herein, as well as all our experiments and results.

And finally, our fourth motivation is a confirming evaluation for the first time of the state-of-the-art Improved NFP-CM model \cite{Rodrigues2017-dy} in the same problem instances reported for NFP-CM \cite[tables 1-2]{Cherri2016-jf} to provide a fair and updated benchmark of the current state-of-the-art MIP models for nesting. We also evaluate the NFP-CMnc model \cite{Cherri2016-jf} in a set of large problem instances not reported before.

\subsection{Definition of the problem and contributions}

The irregular strip-packing problem can be abstractly defined regardless of the geometric representation of the non-overlapping constraints as follows. Let be $\mathcal{P} = \{P_i\}_{i=1,\dots, n}$ the set of pieces to be placed onto the board $\mathcal{B}$, $T = \{t_i\}_{i=1,\dots, n}$ and $R = \{r_i\}_{i=1,\dots, n}$ the set of translation vectors and orientations defining a feasible solution of the problem, $L$ the length of the board, $\mathcal{O}$ the set of admissible orientations, $P(r) \subset \mathbb{R}^2$ denotes an orientation of the piece $P$, and $int(A)$ denotes the interior set for any set $A \subset \mathbb{R}^2$ with the usual topology of $\mathbb{R}^2$. Then, the minimization problem (\ref{model:abstract_problem}) defines the optimal solution to the irregular strip packing problem, where constraints (\ref{ineq:abs-non-overlapping}) prevent the overlapping of pieces, and constraints (\ref{ineq:abs-innerfit}) force the pieces to be entirely contained in the board.
\begin{align}
\text{minimize} \quad & L \label{model:abstract_problem}\\
\text{subject to} \quad & int(P_i(r_i) \oplus t_i) \cap int(P_j(r_j) \oplus t_j) = \varnothing, \quad 1 \leq i < j \leq n  \label{ineq:abs-non-overlapping}\\
& (P_i(r_i) \oplus t_i) \subseteq \mathcal{B}, \quad 1 \leq i \leq n \label{ineq:abs-innerfit}\\
& t_i \in \mathbb{R}^2, r_i \in \mathcal{O}, L \in \mathbb{R}_{>0}
\end{align}

The main research problem tackled in this work is the definition and evaluation of a new and more tightened MIP model for irregular strip packing than the current state-of-the-art family of NFP-CM models \cite{Cherri2016-jf, Rodrigues2017-dy}. Our main contribution is the introduction of a new family of MIP models, called NFP-CM-VS, which is based on a new convex decomposition of the feasible regions between convex parts into vertical slices, together with a new family of valid inequalities, symmetry breakings, and variable eliminations derived from the former geometric decomposition.  Our second significant contribution is the introduction of the first reproducible experimental survey in this line of research, which is based on our software implementation of all models evaluated herein into a Java software library, together with a detailed reproducibility protocol and dataset to allow the exact replication of all our models and results.

The rest of the paper is structured as follows. Section \ref{sec:related_work} reviews the literature on exact mathematical models for irregular strip packing. Section \ref{sec:new_model} introduces our new family of MIP models, whilst the section \ref{sec:evaluatiomn} details our experimental setup and results, and section \ref{sec_discussion} introduces our discussion of the results. Next section summarizes our main conclusions and future work. Subsequently, Appendix A introduces all raw output data generated for each MIP model evaluated herein not included in the results section because of lack of room. Finally, Appendix B introduces a detailed reproducibility protocol based on our supplementary dataset \cite{Lastra-Diaz2022-aw} to allow the exact replication of all our models, experiments, and results. Both aforementioned appendices are provided as supplementary material. 

\section{Related work on exact models for nesting}
\label{sec:related_work}

Our introduction has provided a detailed review of the family of continuous MIP models for irregular strip packing based on the NFP, to which this work belongs. In addition, this section introduces a comprehensive categorization of the literature on exact mathematical models for nesting, together with an extended review of the linear MIP models based on the NFP. However, for a detailed and recent review of the family of mathematical models, we refer the reader to the survey of \citet{Leao2020-bc} and the survey on geometric representations for nesting of \citet{Bennell2008-rc} .

\begin{figure}[t]
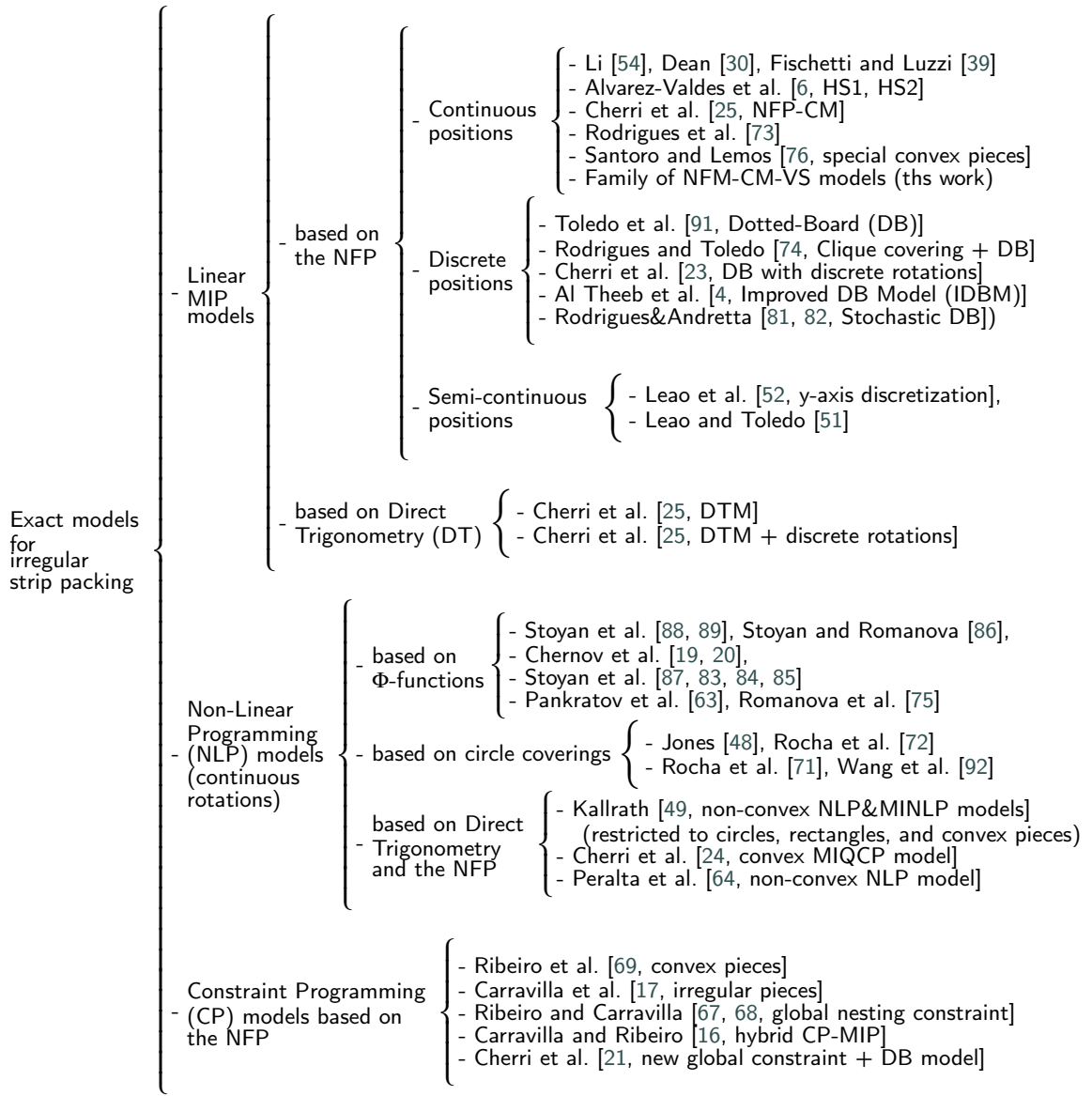

\begin{tabular}{l}
\AB{Exact models \\ for \\ irregular \\ strip packing}
{
    - \AB{Linear \\ MIP \\ models}
    {
        - \AB{based on \\ the NFP}
        {
            - \AB{Continuous \\ positions}
            {
                - \citet{Li1994-co}, \citet{Dean2002-wx}, \citet{Fischetti2009-oa} \\
                - \citet[HS1, HS2]{Alvarez-Valdes2013-wg} \\
                - \citet[NFP-CM]{Cherri2016-jf} \\
                - \citet{Rodrigues2017-dy} \\
                - \citet[special convex pieces]{Santoro2015-px} \\
                - Family of NFM-CM-VS models (ths work)
            }\\
            - \AB{Discrete \\ positions}
            {
                - \citet[Dotted-Board (DB)]{Toledo2013-oi} \\
                - \citet[Clique covering + DB]{Rodrigues2017-zl} \\
                - \citet[DB with discrete rotations]{Cherri2018-rq} \\
                - \citet[Improved DB Model (IDBM)]{Al_Theeb2021-xn} \\
                - Rodrigues\&Andretta \cite[Stochastic DB]{Rodrigues_de_Souza_Queiroz2020-tp, Rodrigues_de_Souza_Queiroz2022-jf}) \\
            }\\
            \bigskip \\
            - \AB{Semi-continuous \\ positions}
            {
                - \citet[y-axis discretization]{Leao2016-xx},\\
                - \citet{Leao2021-zh}
            }
        }\\
        \bigskip \\
        - \AB{based on Direct \\ Trigonometry (DT)}
        {
            - \citet[DTM]{Cherri2016-jf} \\
            - \citet[DTM + discrete rotations]{Cherri2016-jf} \\
        }
    }\\
    \bigskip \\
    - \AB{Non-Linear \\ Programming \\ (NLP) models \\ (continuous \\ rotations)}
    {
        - \AB{based on \\ $\Phi$-functions}
        {
            - \citet{Stoyan1996-kv, Stoyan2012-gu, Stoyan2013-ig}, \\
            - \citet{Chernov2010-ju, Chernov2012-zi}, \\
            - \citet{Stoyan2015-tl, Stoyan2016-ml, Stoyan2016-fg, Stoyan2017-ti} \\
            - \citet{Pankratov2020-ok}, \citet{Romanova2020-sn}
        }\\
        - \AB{based on circle coverings}
        {
            - \citet{Jones2014-fa}, \citet{Rocha2014-vu} \\
            - \citet{Rocha2016-nd}, \citet{Wang2018-yr}
        }\\
        - \AB{based on Direct \\ Trigonometry \\ and the NFP}
        {
            - \citet[non-convex NLP\&MINLP models]{Kallrath2009-jq} \\
            \hspace{0.25cm} (restricted to circles, rectangles, and convex pieces) \\
            - \citet[convex MIQCP model]{Cherri2018-ln}\\
            - \citet[non-convex NLP model]{Peralta2018-ct}
        }\\
    }\\
    \bigskip \\
    - \AB{Constraint Programming \\ (CP) models based on \\ the NFP}
    {
        - \citet[convex pieces]{Ribeiro1999-kg} \\
        - \citet[irregular pieces]{Carravilla2003-hy} \\
        - \citet[global nesting constraint]{Ribeiro2004-nh, Ribeiro2009-gr} \\
        - \citet[hybrid CP-MIP]{Carravilla2005-qb} \\        
        - \citet[new global constraint + DB model]{Cherri2019-nk}
    }
}
\end{tabular}
\caption{Categorization of exact mathematical models for irregular strip packing.}
\label{fig:categorization_models}
\end{figure}

\subsection{Categorization of exact mathematical models}

Figure \ref{fig:categorization_models} shows our categorization of the family of exact continuous mathematical models for irregular strip packing reported in the literature, which can be divided into three large families as follows. First, the family of Linear MIP models based on the NFP, whose main features are the convex decomposition of the feasible regions defined by the pairwise NFP between pieces and the use of mutually exclusive binary variables to select their relative placements, such as the family of continuous models pioneered by \citet[\S8]{Li1994-co}, \citet[\S5]{Dean2002-wx}, and \citet{Fischetti2009-oa}, which are subsequently refined by \citet{Alvarez-Valdes2013-wg}, \citet{Cherri2016-jf}, and \citet{Rodrigues2017-dy}. Second, the family of Constraint Programming (CP) models based on the NFP, whose pioneering work is introduced by \citet{Ribeiro1999-kg}, and subsequently improved by \citet{Carravilla2003-hy}, \citet{Ribeiro2009-gr}, and \citet{Cherri2019-nk}. And third, Other models based on alternatives geometric representations and Non-Linear Programming (NLP) models, such as (3.a) the family of models based on $\Phi$-functions, whose pioneering works are introduced by \citet{Stoyan1996-kv}, \citet{Chernov2010-ju}, and \citet{Stoyan2012-gu}; (3.b) others non-linear models based on direct trigonometry introduced by \citet{Rocha2016-nd}, \citet{Cherri2018-ln}, and \citet{Peralta2018-ct}; and finally, (3.c) the family of models based on circle coverings admitting free rotations, whose pioneering work is introduced by \citet{Jones2014-fa} and subsequently refined by \citet{Rocha2014-vu}, \citet{Rocha2016-nd}, and \citet{Wang2018-yr}. 

Likewise, the linear MIP models can be divided into three subfamilies according to the nature of the decision variables encoding the $(x,y)$ position of the pieces as follows: (1) continuous models, such as the former aforementioned ones \cite{Li1994-co, Dean2002-wx, Fischetti2009-oa,Alvarez-Valdes2013-wg,Cherri2016-jf,Rodrigues2017-dy}; (2) discrete models like the pioneering Dotted-Board model of \citet{Toledo2013-oi}, subsequently refined by \citet{Rodrigues2017-zl}; and (3) semi-continuous models as that proposed by \citet{Leao2016-xx}. On the other hand, the linear MIP models can also be divided into models based on the NFP or direct trigonometry according to the geometric representation used for building the non-overlapping constraints, as shown in figure \ref{fig:categorization_models}.

\subsection{Linear continuous MIP models}

We summarize the review of linear continuous MIP models advanced in our introduction as follows. \citet[\S8]{Li1994-co} introduces the first continuous MIP model for irregular strip packing reported in the literature, whilst \citet[\S5]{Dean2002-wx} evaluates for the first time a minor adaptation of the \citeauthor{Li1994-co}'s model, and \citet{Fischetti2009-oa} introduce and experimentally evaluate a refinement of the \citeauthor{Li1994-co}'s model based on a lifting of the big-M formulation together with a branching-priority algorithm, which is subsequently refined by \citet{Alvarez-Valdes2009-yz} by introducing the HS2 model based on a detailed convex decomposition of the feasible space of relative placements between pieces, a lifting for the bound constraints of the continuous variables, and six new branching strategies. More recently, \citet{Cherri2016-jf} introduce two continuous MIP models, called DTM and NFP-CM, which are based on the convex decomposition of the pieces and the definition of the non-overlapping constraints between pieces by using either direct trigonometry or the convex no-fit polygons among their convex parts, together with some valid inequalities and a variable elimination. Subsequently, \citet{Rodrigues2017-dy} improve the NFP-CM model by breaking the symmetries in the feasible space of relative placements between pieces. On the other hand, \citet{Santoro2015-px} introduce an exact MIP model to compute efficiently tighter upper bounds for the irregular nesting problem based on approximating the input pieces by parallel-chamfered n-gons with up to 8 sides, as shown in \citet[fig. 2]{Santoro2015-px}, and generalizing the disjunctive MIP model for rectangular strip packing proposed by \citet{Sawaya2005-ly}.

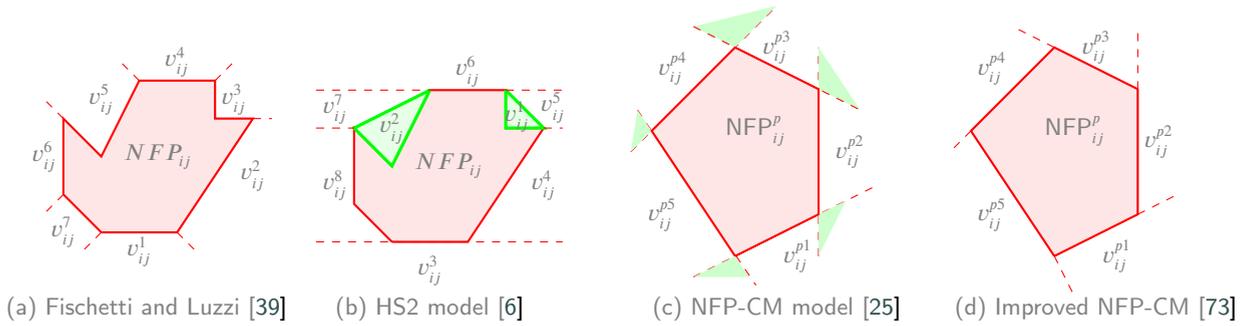
\begin{figure}
\begin{minipage}[t]{0.24\textwidth}
\begin{tikzpicture}[scale=0.5]
\fill[red!10] (1,2) -- (3,2) -- (5,5) -- (4,5) -- (4,6) -- (2,6) -- (1,4) -- (0,5) -- (0,3) -- (1,2);
\draw[thick, red] (1,2) -- (3,2) -- (5,5) -- (4,5) -- (4,6) -- (2,6) -- (1,4) -- (0,5) -- (0,3) -- (1,2);
\draw (2.5,4) node {$NFP_{ij}$};
\draw[dashed, red] (1,2) -- (0.5,1.5);
\draw[dashed, red] (3,2) -- (3.5,1.5);
\draw[dashed, red] (5,5) -- (5.5,5);
\draw[dashed, red] (4,6) -- (4.5,6.5);
\draw[dashed, red] (2,6) -- (1.5,6.5);
\draw[dashed, red] (0,5) -- (-0.5,5.5);
\draw[dashed, red] (0,3) -- (-0.5,2.5);
\draw (2,1.5) node {$v_{ij}^1$};
\draw (2.2,0) node {(a) \citet{Fischetti2009-oa}};
\draw (5,3.5) node {$v_{ij}^2$};
\draw (4.5,5.5) node {$v_{ij}^3$};
\draw (3,6.5) node {$v_{ij}^4$};
\draw (1,5.5) node {$v_{ij}^5$};
\draw (-0.5,4) node {$v_{ij}^6$};
\draw (0,2) node {$v_{ij}^7$};
\end{tikzpicture}
\end{minipage}
\hfill
\begin{minipage}[t]{0.24\textwidth}
\begin{tikzpicture}[scale=0.5]
\fill[red!10] (1,2) -- (3,2) -- (5,5) -- (4,5) -- (4,6) -- (2,6) -- (1,4) -- (0,5) -- (0,3) -- (1,2);
\draw[thick, red] (1,2) -- (3,2) -- (5,5) -- (4,5) -- (4,6) -- (2,6) -- (1,4) -- (0,5) -- (0,3) -- (1,2);
\draw (2.5,4) node {$NFP_{ij}$};
\fill[green!10] (4,6) -- (4,5) -- (5,5) -- (4,6);
\draw[very thick, green] (4,6) -- (4,5) -- (5,5) -- (4,6);
\fill[green!10] (2,6) -- (0,5) -- (1,4) -- (2,6);
\draw[very thick, green] (2,6) -- (0,5) -- (1,4) -- (2,6);
\draw[dashed, red] (-1,6) -- (2,6);
\draw[dashed, red] (4,6) -- (5.5,6);
\draw[dashed, red] (-1,5) -- (0,5);
\draw[dashed, red] (5,5) -- (5.5,5);
\draw[dashed, red] (-1,2) -- (5.5,2);
\draw (4.3,5.3) node {$v_{ij}^1$};
\draw (1,5) node {$v_{ij}^2$};
\draw (2,1.25) node {$v_{ij}^3$};
\draw (2,0.25) node {(b) HS2 model \cite{Alvarez-Valdes2013-wg}};
\draw (5,3.5) node {$v_{ij}^4$};
\draw (5.25,5.55) node {$v_{ij}^5$};
\draw (3,6.5) node {$v_{ij}^6$};
\draw (-0.5,5.5) node {$v_{ij}^7$};
\draw (-0.5,3.5) node {$v_{ij}^8$};
\end{tikzpicture}
\end{minipage}
\hfill
\begin{minipage}[t]{0.24\textwidth}
\begin{tikzpicture}[scale=0.55]
\draw[dashed, red](1.5,5.5) -- (5,9);
\draw[dashed, red](1.66,6.5) -- (4.5,2.3);
\draw[dashed, red](3,2.5) -- (7.33,4.66);
\draw[dashed, red](6,3) -- (6,8);
\fill[green!20] (6,7) -- (7,6.5) -- (6,8) -- (6,7);
\fill[green!20] (4,8) -- (5,9) -- (3,8.5) -- (4,8);
\fill[green!20] (2,6) -- (1.66,6.5) -- (1.5,5.5) -- (2,6);
\fill[green!20] (4,3) -- (3,2.5) -- (4.33,2.5) -- (4,3);
\fill[green!20] (6,4) -- (6,3) -- (6.66,4.33) -- (6,4);
\draw[red,dashed](3,8.5) -- (7,6.5);
\fill[red!10] (4,3) -- (6,4) -- (6,7) -- (4,8) -- (2,6) -- (4,3);
\draw[thick, red] (4,3) -- (6,4) -- (6,7) -- (4,8) -- (2,6) -- (4,3);
\draw (5.5,3) node {$v_{ij}^{p1}$};
\draw (5,1.75) node {(c) NFP-CM model \cite{Cherri2016-jf}};
\draw (6.75,5.5) node {$v_{ij}^{p2}$};
\draw (5,8) node {$v_{ij}^{p3}$};
\draw (2.5,7.5) node {$v_{ij}^{p4}$};
\draw (2.25,4) node {$v_{ij}^{p5}$};
\draw (4.5,6) node {NFP$_{ij}^p$};
\end{tikzpicture}
\end{minipage}
\hfill
\begin{minipage}[t]{0.24\textwidth}
\begin{tikzpicture}[scale=0.55]
\draw[dashed, red](6,4) -- (7, 4.5);
\draw[dashed, red](6,7) -- (6, 8.5);
\draw[dashed, red](4,8) -- (3, 8.5);
\draw[dashed, red](2,6) -- (1.5,5.5);
\draw[dashed, red](4,3) -- (4.5,2);
\fill[red!10] (4,3) -- (6,4) -- (6,7) -- (4,8) -- (2,6) -- (4,3);
\draw[thick, red] (4,3) -- (6,4) -- (6,7) -- (4,8) -- (2,6) -- (4,3);
\draw (5.5,3) node {$v_{ij}^{p1}$};
\draw (5,1.75) node {(d) Improved NFP-CM \cite{Rodrigues2017-dy}};
\draw (6.5,5.75) node {$v_{ij}^{p2}$};
\draw (5,8) node {$v_{ij}^{p3}$};
\draw (2.5,7.5) node {$v_{ij}^{p4}$};
\draw (2.5,4) node {$v_{ij}^{p5}$};
\draw (4.5,6) node {NFP$_{ij}^p$};
\end{tikzpicture}
\end{minipage}
\caption{Disjoint convex decomposition of the feasible relative placements between pieces defined by each model in the family of continuous MIP models for irregular strip packing. Each single convex feasible sub-region in figures (a)-(d) above is enabled by a mutually-exclusive binary variable $v_{ij}^k$. The \citeauthor{Fischetti2009-oa} (a) and HS2 (b) models set their non-overlapping constraints using the overall no-fit polygon $NFP_{ij}$ between pieces $i$ and $j$, whilst the NFP-CP models (c-d) set them using the convex NFP parts between each pair of convex parts of two pieces. Figure (c) shows in green the space of symmetric solutions generated by adjacent feasible sub-regions of the NFP-CM model \cite{Cherri2016-jf}, whilst figure (d) shows the symmetry-breaking of the former model proposed by the Improved NFP-CM model \cite{Rodrigues2017-dy}.}    
\label{fig:convex_decomposition}
\end{figure}

Regardless of the minor differences in the formulation of the continuous MIP models mentioned above, their main difference is the geometric decomposition of the feasible space of relative placements between pieces used to build the non-overlapping constraints, as shown in figure \ref{fig:convex_decomposition}. \citet{Fischetti2009-oa} propose a convex decomposition of the feasible regions using the overall $NFP_{ij}$ between pieces $i$ and $j$, as shown in figure \ref{fig:convex_decomposition}.a. However, \citeauthor{Fischetti2009-oa} do not provide a detailed definition of its geometry, which encourages \citet{Alvarez-Valdes2013-wg} to propose a well-defined convex decomposition based on horizontal slices to build their HS2 model, as shown in figure \ref{fig:convex_decomposition}.b. Subsequently, \citet{Cherri2016-jf} propose a convex decomposition of the pieces and the definition of their non-overlapping constraints by using the convex no-fit polygons among their convex parts, as shown in figure \ref{fig:convex_decomposition}.c. The convex decomposition proposed by \citeauthor{Cherri2016-jf} for their NFP-CM model provides two key advantages on previous MIP models as follows: (1) it avoids the need to explicitly compute the overall NFP between irregular pieces by computing the NFP between convex parts, which is much easier than the former; and, (2) it allows building the non-overlapping constraints between pieces using only one constraint per binary variable, unlike the previous MIP models that require at least three constraints per binary variable.  However, one significant drawback of the NFP-CM model \cite{Cherri2016-jf} is the existence of symmetric solutions induced by the overlapping of adjacent feasible sub-regions, as shown by regions in green in figure \ref{fig:convex_decomposition}.c, which encourages \citet{Rodrigues2017-dy} to break these symmetries by inserting one additional constraint per binary variable activating a separation line defined by the previous edge in the polygon's boundary, as shown in figure \ref{fig:convex_decomposition}.d, at the cost of doubling the number of non-overlapping constraints to break the symmetries in their geometric decomposition.

The main limitation of the family of current continuous MIP models based on the NFP approach is their inability to solve either large problem instances or instances including pieces with complex geometries, which we attribute to two difficulties derived from their structure as follows: (1) the large number of binary variables needed to build the non-overlapping constraints between pieces; and (2) a poor tightening of the LP relaxed  model as a consequence of encoding the linearization of a min-max model in the objective function that does not directly depend on the binary variables. For instance, the number of binary variables in current continuous exact MIP models \cite{Fischetti2009-oa, Alvarez-Valdes2013-wg, Cherri2016-jf, Rodrigues2017-dy} is quadratic regarding the number of pieces $n$ with scalability factor $\mathcal{O}(r^4\frac{n}{2}(n - 1))$, whilst the number of non-overlapping constraints might be up to three times the latter factor, where $r$ is the average number of edges per piece. The $\frac{1}{2}n(n - 1)$ factor is derived from the pairwise combinatory nature of the problem, whilst the $r^4$ gives account of the geometric complexity of the pieces, which derives from the fact that the resulting NFP from two polygons with $s$ and $t$ edges might has $O(s^2t^2)$ edges in the worst case \cite[p.40]{Agarwal2002-ip}, although it is at most $s + t$ for convex polygons. Thus, the complexity of current continuous MIP models grows rapidly regarding both the number of pieces and its geometric complexity.

\subsection{Discrete and semi-continuous MIP models}

For the reasons detailed above, several authors have proposed discrete MIP models to solve efficiently large problem instances as follows: (1) the discrete models based on a grid representation of the board introduced by \citet{Toledo2013-oi}, which is refined by \citet{Rodrigues2017-zl} using a clique-based formulation, being subsequently extended by  \citet{Rodrigues_de_Souza_Queiroz2020-tp, Rodrigues_de_Souza_Queiroz2022-jf} to deal with the uncertainty in the demand of the pieces using a two-stage stochastic programming model; and finally, (2) the semi-continuous model based on the discretization of the y-axis introduced by \citet{Leao2016-xx}. The Dotted-Board model (DB) model of \citet{Toledo2013-oi} is able to solve problem instances with up to 56 pieces of two types to optimality (see \cite[table 3]{Toledo2013-oi}), whilst the semi-continuous MIP model of \citet{Leao2016-xx} removes the resolution error in the x-axis direction and allows solving instances with up to 70 pieces and larger boards \cite[table 1]{Leao2016-xx} than the DB model. However, the discrete models above are approximations of the exact continuous problem. Thus, the quality of their solutions depends on the grid resolution, as shown experimentally by \citet{Sato2016-fd}. Moreover, the computational cost of the discrete models also grows rapidly with the resolution of the grid and the number of different piece types, which encouraged several improvements to the DB model as follows. \citet{Rodrigues2017-zl} propose a new clique-based formulation of the DB model, whilst \citet{Cherri2018-rq} propose two methods to build non-regular grids together with an efficient data structure to represent them. The grid representation, also called no-fit raster, is also used by recent state-of-the-art heuristics methods \cite{Mundim2017-gs, Mundim2018-tk}.

\section{The new family of MIP models for nesting}
\label{sec:new_model}

This section introduces two new MIP models for irregular strip packing based on a new disjoint convex decomposition of the space of feasible placements between pieces defined by the pairwise convex NFP between their convex parts, together with a set of valid inequalities, symmetry breakings, and variable eliminations derived from the former decomposition. The use of the convex NFP between convex parts of the pieces to build the pairwise non-overlapping constraints is the core innovation of the family of NFP-CM models proposed by \citet{Cherri2016-jf} and subsequently improved by \citet{Rodrigues2017-dy} to break the symmetries of the NFP-CM model, as shown in figure \ref{fig:convex_decomposition}.d. Thus, our new MIP models set a more tightened formulation than the NFP-CM model for the same feasible space of relative placements between pieces by defining a new set of constraints and two new binary variables based on our new disjoint convex decomposition shown in figure \ref{fig:vertical_partition}. Our new NFP-CM-VS models include at most 2 binary variables more per convex NFP part than the NFP-CM model \cite{Cherri2016-jf} to encode the left and right feasible sub-regions shown in figure \ref{fig:vertical_partition} in those cases in which the NFP's boundary does not include a vertical edge in some of their sides.

The two aims of our new convex decomposition are as follows: (1) to break the symmetries of the feasible solutions; and (2) to infer a large number of unfeasible relative placements among three pieces or two pieces of the same type, which can be hand-coded either as feasibility cuts or by removing some binary variables from the model, respectively. Our new feasibility cuts among three pieces can also be generalized to combinations of more than three pieces; however, we limit our evaluation herein to piece triplets. First, we introduce the basic formulation of our new MIP model called No-Fit Polygon Covering Model based on Vertical Slices (NFP-CM-VS). Next, we introduce several families of valid inequalities and variable eliminations to break the symmetries and remove some unfeasible solutions derived from the convex NFP parts between two pieces, as well as a family of feasibility cuts among three pieces. Finally, we introduce the NFP-CM-VS2 model integrating all cuts and variable eliminations of the full NFP-CM-VS model but removes the big-M terms for the constraints encoding the vertical slices of our convex decomposition by factorizing all x-axis constraints into two single constraints per convex NFP part.

\subsection{Basic notation and geometric tools}

\emph{Representation of the pieces}. Figure \ref{fig:piece_vars_NFP-CM} details the meaning of the parameters used to represent the geometry of the pieces and the board. The board is represented by an open rectangle of fixed height $H$ and variable length $L$, whose lower and upper bounds are denoted by $L_{lb}$ and $L_{ub}$ respectively. We use a standard orthogonal x-y reference frame setting the origin of the board in its bottom-leftmost corner. Each piece $P_i$ is defined by an irregular polygon decomposed into a set of disjoint convex polygons, called convex parts, whose boundaries are represented by a counter-clockwise ordered list of points in the plane, called \emph{vertices}, such that each pair of consecutive vertices defines a line segment, called \emph{edge}. We use the indexes $f$ and $g$ to denote the source convex parts generating the convex $NFP_{ij}^{fg}$ part between the pieces $P_i$ and $P_j$. For each piece $P_i$, one arbitrary vertex from any of its convex parts is selected as reference point, denoted by $r_i = (x_i,y_i)$, to set the continuous decision variables representing the position of the piece $i$ on the board. The parameters $l_i^{min}$ and $l_i^{max}$ denote the x-axis distance from the leftmost and rightmost points of piece $i$ to its reference point $r_i$, respectively, whilst $h_i^{min}$ and $h_i^{max}$ denote the y-axis distance from the topmost and bottommost points of piece $i$ to $r_i$. The type of each piece $i$ is denoted by $t_i$, whilst their area and demand are denoted by $\Delta_t$ and $d_t$ respectively.

\emph{Convex decomposition of pieces}. To build the non-overlapping constraints of our family of NFP-CM-VS models, all non-convex pieces must be decomposed into convex polygons, as shown in figure \ref{fig:piece_vars_NFP-CM}. Thus, each piece $P_i$ is defined as the union of $Q_i$ convex parts as follows $P_i = \bigcup\limits_{f=1}^{Q_i} P_i^f$. We use the Green's convex decomposition algorithm \cite{Greene1983-fl} implemented by the CGAL library to decompose the pieces into convex parts as proposed by \citet{Cherri2016-jf}, with the aim of replicating their NFP-CM models and experiments. The \citeauthor{Greene1983-fl}'s algorithm efficiently approximates the optimal convex decomposition of a polygon without adding new vertices to its boundary. However, there are many other alternatives in the literature \cite{Fernandez2000-nw} that could potentially generate fewer number of edges, and thus fewer binary variables, than the \citeauthor{Greene1983-fl}'s algorithm. For instance, we found in many problem instances that our implementation of a variant of the Angle-Bisector (AB) convex decomposition method \cite{Agarwal2002-ip} generated models with fewer binary variables than those built using the  \citeauthor{Greene1983-fl}'s method. Thus, any practical implementation of our models should consider both the impact of the convex decomposition methods on the size of the resulting MIP models and their computational cost.

\begin{figure}[t!]
\begin{minipage}[t]{0.45\textwidth}
\begin{tikzpicture}
\fill[black!5] (1.75,2) -- (2,1) -- (3.5,1) -- (4,3) -- (3,4) -- (0.5,3.5) -- (0.5,1.5) -- (1.75,2);
\draw[thick,black] (1.75,2) -- (2,1) -- (3.5,1) -- (4,3) -- (3,4) -- (0.5,3.5) -- (0.5,1.5) -- (1.75,2);
\draw[thick,black, ->] (4,3) -- (3.5,3.5);
\draw[black] (4,3.5) node {$\partial P_i$};
\draw[black] (3,3) node {$P_i^2$};
\draw[black] (2.15,3.35) node {$P_i^1$};
\draw[thin] (1.75,2) -- (3,4);
\draw[dashed,->] (1.75,2) -- (0.5,2);
\draw[dashed,->] (1.75,2) -- (4,2);
\draw[dashed,->] (1.75,2) -- (1.75,4);
\draw[dashed,->] (1.75,2) -- (1.75,1);
\draw[very thick, ->] (0,0.5) -- (6,0.5);
\draw (0,4.5) -- (6,4.5);
\draw[dashed, <->] (5,0.5) -- (5,4.5);
\draw[very thick, ->] (0,0.5) -- (0,4.5);
\draw[dashed, <->] (0,4.1) -- (5,4.1);
\draw (2,4.25) node {$L$};
\filldraw[black] (1.75,2) circle (0.05);
\draw[black] (-0.1,0.15) node {$(0,0)$};
\draw[black] (5.5,0.75) node {$X$};
\draw[black] (-0.25,3.75) node {$Y$};
\draw[black] (5.25,2.5) node {$H$};
\draw[black] (2.6,1.8) node {$r_i=(x_i,y_i)$};
\draw[black] (1.15,2.25) node {$l_i^{min}$};
\draw[black] (3,2.25) node {$l_i^{max}$};
\draw[black] (1.4,3) node {$h_i^{max}$};
\draw[black] (1.4,1.4) node {$h_i^{min}$};
\end{tikzpicture}
\captionof{figure}{Representation of the parameters encoding \\ the geometry of the pieces and the board. Non-convex \\ pieces are decomposed into convex parts, such that \\ the piece $i$ is decomposed as $P_i = P_i^1 \cup P_i^2$.}
\label{fig:piece_vars_NFP-CM}
\end{minipage}
\hfill
\begin{minipage}[t]{0.45\textwidth}
\begin{tikzpicture}[scale=0.75]
\fill[black!5] (4,3) -- (6,3) -- (8,6) -- (8,8) -- (2,8) -- (2,6) -- (4,3);
\draw[thick, black] (4,3) -- (6,3) -- (8,6) -- (8,8) -- (2,8) -- (2,6) -- (4,3);
\draw[thick, black, ->] (2,6) -- (3,4.5);
\draw[black] (2,4.5) node {$\partial NFP_{AB}$};
\fill[black!15] (6,5) -- (8,8) -- (4,8) -- (6,5);
\draw[black] (6,5) -- (8,8) -- (4,8) -- (6,5);
\fill[black!15] (8.5,4) -- (10.5,4) -- (10.5,6) -- (8.5,6) -- (8.5,4);
\draw[black] (8.5,4) -- (10.5,4) -- (10.5,6) -- (8.5,6) -- (8.5,4);
\filldraw[black] (6,5) circle (0.1);
\filldraw[black] (8.5,4) circle (0.1);
\draw[black] (6,7) node {$A$};
\draw[black] (9.5,5) node {$B$};
\draw[black] (4,5.5) node {$NFP_{AB}$};
\draw[black] (6,4.7) node {$r_A$};
\draw[black] (8.5,3.7) node {$r_B$};
\draw[black] (4.5,6.5) node {$\partial A$};
\draw[black, ->] (4,8) -- (5,6.5);
\draw[black] (9.5,6.35) node {$\partial B$};
\draw[black,->] (10.5,6) -- (9.5,6);
\end{tikzpicture}
\captionof{figure}{No-fit polygon $NFP_{AB}$ between the static \\ polygon $A$ and the orbiting polygon $B$, whose boundaries \\are denoted by $\partial NFP_{AB}$, $\partial A$, and $\partial B$, respectively.} 
\label{fig:nofit_polygon}
\end{minipage}
\\
\begin{minipage}[t]{0.45\textwidth}
\begin{tikzpicture}
\draw[thick, black] (-1,6) -- (7,6);
\draw[thick, black] (-1,0) -- (7,0);
\draw[thick, black] (-1,0) -- (-1,6);
\draw[thick, black] (7,0) -- (7,6);
\fill[black!5] (-1,0) -- (-1,6) -- (7,6) -- (7,0) -- (-1,0);
\fill[black!15] (2,1) -- (3,1) -- (4,3) -- (3,4.5) -- (2,4) -- (1,3) -- (1,2) -- (2,1);
\draw[black] (2,1) -- (3,1) -- (4,3) -- (3,4.5) -- (2,4) -- (1,3) -- (1,2) -- (2,1);
\draw[dashed, <->] (-1,5.5) -- (7,5.5);
\draw (0.25, 5.75) node {$L$};
\draw[dashed, <->] (-0.5,0) -- (-0.5,6);
\draw (-0.25, 5) node {$H$};
\draw[thin, black] (1,0) -- (1,6);
\draw[thin, black] (4,0) -- (4,6);
\draw[thin, black] (3,4.5) -- (3,6);
\draw[thin, black] (2,4) -- (2,6);
\draw[thin, black] (2,1) -- (2,0);
\draw[thin, black] (3,1) -- (3,0);
\draw[black] (0.25, 3) node {$v_{ij}^{fgl}$};
\draw[black] (5.5, 3) node {$v_{ij}^{fgr}$};
\draw[black] (2.5, 2.5) node {$NFP_{ij}^{fg}$};
\draw[black] (1.5,0.5) node {$v_{ij}^{fgb_1}$};
\draw[black] (2.5,0.5) node {$v_{ij}^{fgb_2}$};
\draw[black] (3.5,0.5) node {$v_{ij}^{fgb_3}$};
\draw[black] (1.5,4.5) node {$v_{ij}^{fgt_3}$};
\draw[black] (2.5,5) node {$v_{ij}^{fgt_2}$};
\draw[black] (3.6,4.75) node {$v_{ij}^{fgt_1}$};
\filldraw[black] (2,1) circle (0.05);
\filldraw[black] (3,1) circle (0.05);
\filldraw[black] (4,3) circle (0.05);
\filldraw[black] (3,4.5) circle (0.05);
\filldraw[black] (2,4) circle (0.05);
\filldraw[black] (1,3) circle (0.05);
\filldraw[black] (1,2) circle (0.05);
\end{tikzpicture}
\caption{Disjoint convex decomposition into vertical slices \\ of the feasible region for the relative placement of piece $P_j$ \\ regarding piece $P_i$, as defined by the convex no-fit polygon \\ $NFP_{ij}^{fg}$ in dark grey, which are enabled in our MIP models \\ by the set of mutually-exclusive binary variables $v_{ij}^{fgk}$.}
\label{fig:vertical_partition}
\end{minipage}
\hfill
\begin{minipage}[t]{0.45\textwidth}
\begin{tikzpicture}[scale = 1.30]
\fill[black!5] (3,1.5) -- (4,1.5) -- (4.5,2) -- (4.5,3) -- (3,3) -- (2.5,2.5) -- (2.5,2) -- (3,1.5);
\draw[black] (3,1.5) -- (4,1.5) -- (4.5,2) -- (4.5,3) -- (3,3) -- (2.5,2.5) -- (2.5,2) -- (3,1.5);
\fill[black!5] (4.5,3) -- (4.5,4) -- (4,4.5) -- (3,4.5) -- (2.5,4) -- (2.5,3.5) -- (3,3);
\draw[black] (4.5,3) -- (4.5,4) -- (4,4.5) -- (3,4.5) -- (2.5,4) -- (2.5,3.5) -- (3,3);
\draw[dashed, black] (3,3) -- (4.5,3);
\filldraw[black] (3,1.5) circle (0.05);
\draw[black]  (3.5,3.75) node {$A^1$};
\draw[black]  (3.5,2.25) node {$A^2$};
\fill[black!5] (5.5,2.5) -- (6.5,2.5) -- (6.5,3.5) -- (5.5,3.5) -- (5.5,2.5);
\draw[black] (5.5,2.5) -- (6.5,2.5) -- (6.5,3.5) -- (5.5,3.5) -- (5.5,2.5);
\filldraw[black] (5.5,2.5) circle (0.05);
\draw[black]  (6,3) node {$B^1$};
\draw[thick, dashed, black] (2,0.5) -- (4,0.5) -- (4.5, 1.0) -- (4.5,3) -- (2,3) -- (1.5,2.5) -- (1.5,1) -- (2,0.5);
\draw[thick, dashed, black] (2,2) -- (4.5,2) -- (4.5,4) -- (4,4.5) -- (2,4.5) -- (1.5,4) -- (1.5,2.5) -- (2,2);
\draw[black] (3,1) node {$NFP_{AB}^{21}$};
\draw[black]  (2.05,3.25) node {$NFP_{AB}^{11}$};
\draw[black]  (6,2) node {$A = A^1 \cup A^2$};
\draw[black]  (5.75,1.5) node {$B = B^1$};
\end{tikzpicture}
\caption{Dashed lines show the boundaries of the \\ convex $NFP_{AB}^{11}$ and $NFP_{AB}^{21}$ parts between the single \\ convex part $B^1$ of piece $B$ and the two convex parts $A^1$ \\ and $A^2$ of piece $A$.}
\label{fig:convex_NFP_parts}
\end{minipage}
\end{figure}

\emph{Computation of convex NFP parts}. Figure \ref{fig:nofit_polygon} shows in light grey the no-fit polygon $NFP_{AB}$ between two convex polygons $A$ and $B$ with reference points $r_A$ and $r_B$, respectively. The outer and boundary regions of $NFP_{AB}$ define the feasible region in which polygon $B$ can be placed without overlapping polygon $A$, whilst its inner region sets the non-feasible relative positions for placing $B$, and the boundary of $NFP_{AB}$, denoted by $\partial NFP_{AB}$, sets the positions in which both polygons are in contact. Given two convex polygons $A$ and $B$, their no-fit polygon $NFP_{AB} = A \oplus (-B(0,0))$ is always convex \cite[theorem 13.5]{De_Berg1997-bi} and it can be efficiently computed using any specialized algorithm for convex polygons, such as the orbiting method of \citet{Cuninghamegreen1989-si}, or a specialized version of Minkowski sums for convex polygons \cite[p.299]{De_Berg1997-bi}. Despite  \citet{Cherri2016-jf} propose the use of the \citeauthor{Cuninghamegreen1989-si}'s algorithm to compute the NFP between convex polygons to build their NFP-CM model, our preferred option is to use convex Minkowski sums because once all convex parts of any orbiting piece $B$ are translated to the origin by summing the $-r_{B}$ vector, the relative positions of all convex $NFP_{AB}^{fg}$ parts are well defined regarding the reference point $r_{A}$ of the static piece, unlike the resulting NFP parts obtained with the \citeauthor{Cuninghamegreen1989-si}'s algorithm. Algorithm \ref{alg:Minkowski} introduces our detailed implementation of the \emph{MinkowskiSum} algorithm \cite[p.299]{De_Berg1997-bi} for convex polygons used to compute the convex NFP between the convex parts of each pair of pieces.

\begin{algorithm}[t!]
\caption{Our version of the MinkowskiSum$(\delta_a,\delta_b)$ algorithm for convex polygons introduced by \citet[p.299]{De_Berg1997-bi}. ConvexMinkowskiSum function is used to compute the convex no-fit polygon $NFP_{AB} = A \oplus (-B(0,0))$ between convex parts $A$ and $B$ by calling the function below with parameters $\delta_a =\partial A$ and $\delta_b=-(\partial B \oplus -r_B)$.}
\label{alg:Minkowski}
\begin{algorithmic}[1]
\Require counter-clockwise oriented boundaries $\delta_a \in \mathbb{R}^{2 \times n}, \delta_b \in \mathbb{R}^{2 \times m}$
\Ensure counter-clockwise boundary $\delta_{NFP}$ of $NFP_{AB}$
\Function{ConvexMinkowskiSum}{$\delta_a, \delta_b$}
\State $\delta'_a \gets sort(\delta_a)$
\State $\delta'_b \gets sort(\delta_b)$ \Comment{\text{first vertex has the smallest y-axis coordinate}}
\State n $\gets$ lengthof$(\delta'_a)$
\State m $\gets$ lengthof$(\delta'_b)$
\State $\delta_{NFP} \gets \emptyset$ \Comment{$\delta_{NFP}$ \text{ is an array of } $\mathbb{R}^{2}$ \text{points}}
\State i $\gets$ 0
\State j $\gets$ 0
\While{$(i < n \lor j < m)$}
\State $\delta_{NFP} \: \xleftarrow{adds} \: \delta'_a[i\%n] + \delta'_b[j\%m]$ \Comment{\text{adds a new vertex to $\delta_{NFP}$}}
\State $b \gets \delta'_a[(i+1)\%n] - \delta'_a[i\%n]$
\State $c \gets \delta'_b[(j+1)\%m] - \delta'_b[j\%m]$
\State $\theta \gets b_x c_y - c_x b_y$
\If {$\theta \geq 0$}
\State $i \gets i + 1$
\EndIf
\If {$\theta \leq 0$}
\State $j \gets j + 1$
\EndIf
\EndWhile
\State \Return $\delta_{NFP}$
\EndFunction
\end{algorithmic}
\end{algorithm}

\emph{New convex decomposition based on vertical slices}. Figure \ref{fig:vertical_partition} shows our new disjoint convex decomposition of the feasible space of relative placements for the piece $j$ regarding the piece $i$, such that their corresponding convex parts $g$ and $f$ do not overlap. We recall that the non-overlapping constraints of our family of NFP-CM-VS model are defined by all pairwise combination of convex parts between two distinct pieces. $NFP_{ij}^{fg}$ denotes the convex no-fit polygon between the static convex part $f$ of piece $i$ and the orbiting convex part $g$ of piece $j$. For example, figure \ref{fig:convex_NFP_parts} shows the boundaries of $NFP_{AB}^{11}$ and $NFP_{AB}^{21}$ parts derived from the combinations of convex parts of pieces $A$ and $B$ as dashed lines, whilst $\partial NFP_{AB}^{11}$ and $\partial NFP_{AB}^{21}$ denote their boundaries. Likewise, $\partial NFP_{ij}^{fg}$ denotes the boundary of $NFP_{ij}^{fg}$ and it is defined by a counter-clockwise oriented polyline whose line segments are called \emph{edges} and denoted by $(a_{ij}^{fgk}, b_{ij}^{fgk}) \in \mathbb{R}^2 \times \mathbb{R}^2$, with extreme points $a_{ij}^{fgk} = (a_{ij,x}^{fgk}, a_{ij,y}^{fgk})$ and $b_{ij}^{fgk} = (b_{ij,x}^{fgk}, b_{ij,y}^{fgk})$.

The no-fit boundary $\partial NFP_{ij}^{fg}$ is represented as the union of all its line segments, called edges, as defined in expression (\ref{def:nfp_boundary}) below, being $\mathcal{K}_{ij}^{fg}$ the overall number of boundary edges. To define the feasible sub-regions induced by $NFP_{ij}^{fg}$, we decompose the line segments $e_{ij}^{fgk}$ of $\partial NFP_{ij}^{fg}$ into three disjoint sets grouping the top, bottom, and side edges, as defined by the expressions (\ref{def:nfp_boundary_decomposition}-\ref{set:side_edges}) below.
\begin{align}
\partial NFP_{ij}^{fg} &= \bigcup\limits_{k=1}^{\mathcal{K}_{ij}^{fg}} e_{ij}^{fgk}, \quad e_{ij}^{fgk} = (a_{ij}^{fgk}, b_{ij}^{fgk}), \quad a_{ij}^{fgk}, b_{ij}^{fgk} \in \mathbb{R}^2 \label{def:nfp_boundary} \\
\partial NFP_{ij}^{fg} &= \mathcal{T}_{ij}^{fg} \cup \mathcal{B}_{ij}^{fg} \cup \mathcal{S}_{ij}^{fg} \qquad\qquad\qquad\qquad \text{(boundary's edge decomposition)} \label{def:nfp_boundary_decomposition} \\
\mathcal{T}_{ij}^{fg} &= \{e_{ij}^{fgk} \in \partial NFP_{ij}^{fg} | a_{ij,x}^{fgk} > b_{ij,x}^{fgk}\} \qquad \text{(top edges)} \label{set:top_edges} \\
\mathcal{B}_{ij}^{fg} &= \{e_{ij}^{fgk} \in \partial NFP_{ij}^{fg} | a_{ij,x}^{fgk} < b_{ij,x}^{fgk}\} \qquad \text{(bottom edges)} \label{set:bottom_edges} \\
\mathcal{S}_{ij}^{fg} &= \{e_{ij}^{fgk} \in \partial NFP_{ij}^{fg} | a_{ij,x}^{fgk} = b_{ij,x}^{fgk}\} \qquad \text{(side edges)} \label{set:side_edges}
\end{align}

\newpage
Binary variables enabling top feasible regions defined by edges in $\mathcal{T}_{ij}^{fg}$ are denoted by $v_{ij}^{fgt_k}$ in figure \ref{fig:vertical_partition}, whilst binary variables enabling bottom feasible regions defined by edges in $\mathcal{B}_{ij}^{fg}$ are denoted by $v_{ij}^{fgb_k}$. Finally, side edges $\mathcal{S}_{ij}^{fg}$ will not be considered in our model, because the feasible regions defined by vertical edges, either on the left or right sides of $NFP_{ij}^{fg}$, are enabled by the distinguished binary variables $v_{ij}^{fgl}$ and $v_{ij}^{fgr}$, as shown in figure \ref{fig:vertical_partition}. Each convex feasible sub-region enabled by a variable $v_{ij}^{fgk}$ is denoted by $R_{ij}^{fgk} \subset \mathbb{R}^2$. Finally, we introduce a notation for the smallest and largest x-axis coordinates of each convex $NFP_{ij}^{fg}$ part, as detailed in equations (\ref{eq:xmin_ij}) and (\ref{eq:xmax_ij}) below.
\begin{align}
\underline{x}_{ij}^{fg} &= \min\{\partial NFP_{ij,x}^{fg}\} \qquad \text{(smallest x-axis coordinate of the convex NFP part)} \label{eq:xmin_ij} \\
\overline{x}_{ij}^{fg} &= \max\{\partial NFP_{ij,x}^{fg}\} \qquad \text{(largest x-axis coordinate of the convex NFP part)} \label{eq:xmax_ij}
\end{align}

\emph{Definition of index sets}. Before detailing our models, we define the index set $I_{ij}^{fg}$ detailed below (\ref{eq:index_set}) to simplify our notation and the presentation of our models. Each tuple $(i,j,f,g) \in I_{ij}^{fg}$ is used to denote the indexes involved in the definition of the binary variables or constraints concerning a convex $NFP_{ij}^{fg}$ part obtained by evaluating the algorithm \ref{alg:Minkowski} with the convex parts $f$ and $g$ of pieces $i$ and $j$ as input, where $Q_i$ and $Q_j$ are the number of convex parts of the former pieces, respectively. On the other hand, $T_{ij}^{fg}$ and $B_{ij}^{fg}$ denote the index set of the binary variables encoding the top and bottom convex feasible sub-regions defined by $NFP_{ij}^{fg}$, whilst $K_{ij}^{fg}$ denotes the index set of all binary variables encoding the convex feasible sub-regions shown in figure \ref{fig:vertical_partition}, as detailed below.
\begin{align}
I_{ij}^{fg} &= \{\alpha \in \{1,\dots,N\} \times \{1,\dots,N\} \times \{1,\dots,Q_i\} \times \{1,\dots,Q_j\} \:|\: 1 \leq i < j \leq N \} \label{eq:index_set} \\
T_{ij}^{fg} &= \{k \in \{1,\dots,\mathcal{K}_{ij}^{fg}\} \:|\: e_{ij}^{fgk} \in \mathcal{T}_{ij}^{fg}\} \: \text{(indexes of variables encoding top feasible sub-regions)} \label{lab:top_indexes} \\
B_{ij}^{fg} &= \{k \in \{1,\dots,\mathcal{K}_{ij}^{fg}\} \:|\: e_{ij}^{fgk} \in \mathcal{B}_{ij}^{fg}\} \: \text{(indexes of variables encoding bottom feasible sub-regions)} \label{lab:bottom_indexes} \\
K_{ij}^{fg} &= \{k \in T_{ij}^{fg} \cup B_{ij}^{fg} \cup \{v_{ij}^{fgl}, v_{ij}^{fgr}\}\} \: \text{(indexes of variables encoding all feasible sub-regions)} \label{lab:all_indexes}
\end{align}

Most of the constraints of our models are based on one convex NFP part resulting from the combination of two convex parts from two different pieces, whose definition requires a tuple with five indexes from the Cartesian product of $I_{ij}^{fg}$ with any of the index sets of binary variables in definitions (\ref{lab:top_indexes}-\ref{lab:all_indexes}). However, we also introduce here for the first time several constraints among three pieces that demand up to nine indexes for their definition. Thus, we define in (\ref{eq:index_set2}) the index set $II_{iju}^{fgh}$ for the convex parts $f,g$, and $h$ from three different pieces denoted by $i,j$, and $u$ respectively.
\begin{align}
II_{iju}^{fgh} &= \{\alpha \in \{1,\dots,N\} \times \{1,\dots,N\} \times \{1,\dots,N\} \times \{1,\dots,Q_i\} \times \{1,\dots,Q_j\} \times \{1,\dots,Q_u\} \nonumber \\
&  \:|\: 1 \leq i < j < u\leq N \} \label{eq:index_set2}
\end{align}

\subsection{The NFP-CM-VS models}

The basic NFP-CM-VS model is defined by the objective function (\ref{model:nfp_cm_vs_objective}) and the constraints (\ref{model:x_bounds}-\ref{var:domain_binary_vars}), whilst the full NFP-CM-VS model is defined by the former objective function and the constraints, symmetry breakings, valid inequalities, and variable eliminations in expressions (\ref{model:x_bounds}-\ref{ineq:unfeasibility_cuts_for_triplets}). The objective function (\ref{model:nfp_cm_vs_objective}) together with the constraints (\ref{model:x_bounds}-\ref{eq:binary_variables}) fit the definition of the basic NFP-CM model without cuts \cite{Cherri2016-jf}, with the only exception of our two distinguished binary variables $v_{ij}^{fgl}$ and $v_{ij}^{fgl}$ encoding the left and right feasible sub-regions shown in figure \ref{fig:vertical_partition}, whilst the constraints (\ref{model:nfp_cm_vs_first_constraint}-\ref{model:nfp_cm_vs_last_constraint}) encode our new convex decomposition based on vertical slices, which removes all symmetric solutions derived from any relative placement between pieces. $N$ and $m$ denote the number of pieces and types of pieces, respectively.

Constraints (\ref{model:x_bounds}) and (\ref{model:y_bounds}) set the lower and upper bounds for the reference point $r_i = (x_i, y_i)$ of each piece to ensure that all pieces are inside the board. These two later constraints encode the \emph{Inner-Fit Polygon} (IFT) of each piece, which represents the feasible region of the board in which it can be placed. Constraint (\ref{ineq:L_upperbound}) sets the upper bound $L_{ub}$ for $L$ defined as the overall sum of the length of all pieces, whilst constraint (\ref{ineq:L_lowerbound}) sets the lower bound $L_{lb}$ of $L$, which is defined as the largest value between the largest length of any piece and the overall sum of the areas of the pieces divided by the height $(H)$ of the board. The lower and upper bounds of $L$ detailed above have been also used by most of continuous MIP models reported in the literature \cite{Fischetti2009-oa, Alvarez-Valdes2013-wg, Cherri2016-jf, Rodrigues2017-dy}. Despite having been studied other upper and lower bounds for $L$, the bounds mentioned above have become the standard ones for the objective function of the irregular strip-packing problem because of the drawbacks of other alternatives. For instance, \citet{Alvarez-Valdes2013-wg} propose a lower bound for $L$ based on the approximation of the pieces by a set of inner rectangles and the solution of an integer program defining a 1-Contiguous Bin Packing Problem (1-CBPP) \cite{Alvarez-Valdes2009-yz}. However, the authors finally discard the use of this later lower bound because its computational cost is very high. On the other hand, \citet{Cherri2016-jf} point out on the possibility of using more tightened upper bounds that ``a tighter big-M formulation generally makes it hard to the solver to find good quality solutions at the beginning of the search", a conclusion that we also subscribe. Thus, this later drawback and the computational cost of computing tighter upper bounds endorse using the upper bound introduced above as a practical solution.
\begin{align}
\text{min} \quad & L \label{model:nfp_cm_vs_objective} \\
\text{s.t.} \quad & l_i^{min} \leq x_i \leq L - l_i^{max} \quad\quad\quad\quad 1 \leq i \leq N \label{model:x_bounds} \\
 & h_i^{min} \leq y_i \leq H - h_i^{max} \quad\quad\quad 1 \leq i \leq N \label{model:y_bounds} \\
 & L \leq \sum_{t=1}^m d_t(l_t^{min} + l_t^{max}) = L_{ub} \label{ineq:L_upperbound} \\ 
 & L_{lb} = \max\{\max\limits_{1 \leq t \leq m}\{l_t^{min} + l_t^{max}\}, \frac{1}{H} \sum_{t=1}^m d_t \Delta_t\} \leq L \label{ineq:L_lowerbound} \\
 & (b_{ij,x}^{fgk} - a_{ij,x}^{fgk})(y_j - y_i) + (a_{ij,y}^{fgk} - b_{ij,y}^{fgk})(x_j - x_i)  + C_{ij}^{fgk} \leq (1 - v_{ij}^{fgk})M_{ij}^{fgk} \nonumber \\
 & \qquad\qquad   \forall (i,j,f,g,k) \in \{I_{ij}^{fg} \times T_{ij}^{fg} \cup B_{ij}^{fg}\} \label{ineq:main_edge_feasible_region}  \\
 & \quad\quad\quad\quad C_{ij}^{fgk} = b_{ij,y}^{fgk}a_{ij,x}^{fgk} - b_{ij,x}^{fgk}a_{ij,y}^{fgk} \nonumber \\
  & \quad\quad\quad\quad M_{ij}^{fgk} \geq |b_{ij,x}^{fgk} - a_{ij,x}^{fgk}|H + |a_{ij,y}^{fgk} - b_{ij,y}^{fgk}|L_{ub} + C_{ij}^{fgk}  \nonumber \\
 & v_{ij}^{fgl} + v_{ij}^{fgr} + \sum_{k \in T_{ij}^{fg} \cup B_{ij}^{fg}} v_{ij}^{fgk} = 1, \qquad \forall (i,j,f,g) \in I_{ij}^{fg} \label{eq:binary_variables} \\
 & b_{ij,x}^{fgk} + x_i - x_j \leq (1 - v_{ij}^{fgk}) M_{ij}^{\prime fgk}, \qquad \forall (i,j,f,g,k) \in \{I_{ij}^{fg} \times T_{ij}^{fg}\}  \label{model:nfp_cm_vs_first_constraint} \\
 & \qquad M_{ij}^{\prime fgk} \geq b_{ij,x}^{fgk} + L_{ub} - l_i^{max} - l_j^{min} \nonumber \\
 & x_j - x_i - a_{ij,x}^{fgk} \leq (1 - v_{ij}^{fgk}) M_{ij}^{\prime\prime fgk}, \qquad \forall (i,j,f,g,k) \in \{I_{ij}^{fg} \times T_{ij}^{fg}\} \\
 & \qquad M_{ij}^{\prime\prime fgk} \geq L_{ub} - l_j^{max} - l_i^{min} - a_{ij,x}^{fgk} \nonumber \\
 & a_{ij,x}^{fgk} + x_i - x_j \leq (1 - v_{ij}^{fgk}) \bar{M}_{ij}^{\prime fgk}, \qquad \forall (i,j,f,g,k) \in \{I_{ij}^{fg} \times B_{ij}^{fg}\} \\
 & \qquad \bar{M}_{ij}^{\prime fgk} \geq a_{ij,x}^{fgk} + L_{ub} - l_i^{max} - l_j^{min} \nonumber \\
 & x_j - x_i - b_{ij,x}^{fgk} \leq (1 - v_{ij}^{fgk}) \bar{M}_{ij}^{\prime\prime fgk}, \qquad \forall (i,j,f,g,k) \in \{I_{ij}^{fg} \times B_{ij}^{fg}\} \\
 & \qquad \bar{M}_{ij}^{\prime\prime fgk} \geq L_{ub} - l_j^{max} - l_i^{min} - b_{ij,x}^{fgk} \nonumber \\
 & x_j - x_i - \underline{x}_{ij}^{fg} \leq (1 - v_{ij}^{fgl})M_{ij}^{fgl}, \qquad \forall (i,j,f,g) \in I_{ij}^{fg} \\
 & \qquad M_{ij}^{fgl} \geq L_{ub} - l_j^{max} - l_i^{min} - \underline{x}_{ij}^{fg} \nonumber \\
 & x_j - x_i - \overline{x}_{ij}^{fg} \geq (1 - v_{ij}^{fgr})M_{ij}^{fgr}, \qquad \forall (i,j,f,g) \in I_{ij}^{fg} \label{model:nfp_cm_vs_last_constraint}\\
 & \qquad M_{ij}^{fgr} \leq l_j^{min} + l_i^{max} - L_{ub} - \overline{x}_{ij}^{fg} \nonumber \\
 & \qquad  L \in \mathbb{R}_{>0}  \label{var:domain_L} \\
 & \qquad (x_i,y_i) \in \mathbb{R}^2,  \; 1 \leq i \leq N  \label{var:domain_xy} \\
 & \qquad v_{ij}^{fgk} \in \{0,1\}, \qquad \forall (i,j,f,g,k) \in \{I_{ij}^{fg} \times K_{ij}^{fg}\}
\label{var:domain_binary_vars}
\end{align}

\newpage
Constraints (\ref{ineq:main_edge_feasible_region}) defines the feasible region for any feasible relative placement of piece $j$ regarding piece $i$, which is defined on the right side of each counter-clockwise oriented edge $e_{ij}^{fgk} \in \partial NFP_{ij}^{fg}$ for all convex NFP parts between both former pieces. The definition of the constraints (\ref{ineq:main_edge_feasible_region}) on the convex $NFP_{IJ}^{fg}$ parts between pieces is the core contribution of the NFP-CM model introduced by \citet{Cherri2016-jf}, which has the advantage of requiring the minimum number of constraints to define the feasible space of relative placements between pieces among the current family of continuous MIP models. However, the feasible space generated by the constraints (\ref{ineq:main_edge_feasible_region}) introduces a large number of symmetric solutions derived from the overlapping of feasible regions spanned by consecutive edges, as shown in figure \ref{fig:convex_decomposition}.c. In order to break these later symmetries, \citet{Rodrigues2017-dy} propose the convex decomposition of the feasible space shown in figure \ref{fig:convex_decomposition}.c, which doubles the number of non-overlapping constraints per binary variable of the NFP-CM model \cite{Cherri2016-jf}. Unlike the Improved NFP-CM model of \citet{Rodrigues2017-dy}, our novel disjoint convex decomposition shown in figure \ref{fig:vertical_partition} allows the NFP-CM-VS model to break the symmetries of the solution space at the expense of tripling the number of constraints per binary variable of the NFP-CM model \cite{Cherri2016-jf}. However, our new formulation allows the explicit derivation of many feasibility cuts and symmetry breaks linking the binary variables encoding the relative placements among multiple pieces, which provide a tighter formulation of the problem than the former state-of-the-art NFP-CM models \cite{Cherri2016-jf, Rodrigues2017-dy}. On the other hand, the formulation of our NFP-CM-VS2 model introduced in section \ref{sec:NFP_CM_VS2} removes the aforementioned drawback on the number of constraints required by the symmetry-breaking of our NFP-CM-VS model.

Constraints (\ref{eq:binary_variables}) encode the selection in any feasible solution of a single feasible sub-region for the relative placement between two pieces, as defined by our convex decomposition of each convex $NFP_{ij}^{fg}$ part shown in figure \ref{fig:vertical_partition}. Constraints (\ref{model:nfp_cm_vs_first_constraint}-\ref{model:nfp_cm_vs_last_constraint}) encode the vertical lines delimiting our convex feasible sub-regions, as shown in figure \ref{fig:vertical_partition}. Our model defines the same binary variables for the edges of $\partial NFP_{ij}^{fg}$ than the NFP-CM \cite{Cherri2016-jf} and Improved NFP-CM \cite{Rodrigues2017-dy} models, with the only exception of the two distinguished binary variables $v_{ij}^{fgl}$ and $v_{ij}^{fgr}$ encoding the left and right sub-regions, as shown in figure \ref{fig:vertical_partition}. Finally, constraints (\ref{var:domain_L} - \ref{var:domain_binary_vars}) set the domains for the decision variables.

\subsubsection{Symmetry-breaking and variable eliminations for identical pieces}

\emph{Symmetry-breaking of piece permutations}. It is a well-known fact that any permutation of identical pieces generates the same overall feasible space and optimal solutions, which leads the Branch and Bound (B\&B) algorithm to obtain many symmetrical solutions. To break the aforementioned symmetries, we can impose that either $x_i \leq x_j$ for all $1 \leq i < j \leq N$ if pieces $i$ and $j$ are of the same type, as proposed by \citet[\S4.5]{Alvarez-Valdes2013-wg} and \citet[\S3.1]{Cherri2016-jf}, or $y_i \leq y_j$ as proposed by \citet[ineq.15]{Rodrigues2017-dy}. Our family of NFP-CM-VS models includes the y-axis symmetry-breaking proposed for the Improved NFP-CM model \cite{Rodrigues2017-dy}, as defined by constraints (\ref{ineq:yaxis_symmetry}) below. However, the satisfaction of the constraints (\ref{ineq:yaxis_symmetry}) joined to our new convex decomposition shown in figure \ref{fig:vertical_partition} induces several a priori unfeasible relative placements between identical pieces that can be coded into the models as a set of new variable eliminations and valid inequalities, as detailed by constraints (\ref{eq:zero_variables}) and (\ref{ineq:new_inequality_identical_pieces}) below.
\begin{equation}
y_i \leq y_j,\quad 1 \leq i < j \leq N, j = \min \{k = i + 1,\dots, N \: | \: t_i = t_j \}
\label{ineq:yaxis_symmetry}
\end{equation}

\emph{Variable eliminations for identical pieces.} Because constraint (\ref{ineq:yaxis_symmetry}) must be satisfied for identical pieces, all binary variables $v_{ij}^{fgk}$ enabling the bottom regions of any $NFP_{ij}^{fg}$ between two identical pieces must be equal to 0 whenever pieces $i$ and $j$ are of the same type (i.e. $t_i = t_j$) and the bottom region is below the reference point $r_i$. For this reason, all binary variables $v_{ij}^{fgk}$ encoding bottom feasible sub-regions should be fixed to 0 and removed from any NFP-CM-VS model, as detailed by constraints (\ref{eq:zero_variables}) below. Because the origin $(0,0)$ of all convex $NFP_{ij}^{fg}$ parts is defined at the reference point $r_i$, a bottom feasible region $R_{ij}^{fgk}$ will be below $r_i$ if the y-coordinates of the extreme points of its associate edge $e_{ij}^{fgk}$ are negative, as detailed by equality (\ref{eq:zero_variables}).

\begin{equation}
    v_{ij}^{fgk} = 0  \qquad \forall (i,j,f,g,k) \in \{I_{ij}^{fg} \times B_{ij}^{fg} \:|\: t_i = t_j \land (a_{if,y}^{fg} < 0) \land (b_{ij,y}^{fg} < 0)\} \label{eq:zero_variables}
\end{equation} 

\emph{Valid inequality for identical pieces}. For the same reasons detailed above, given two orbiting identical pieces $j$ and $u$ and another static distinct piece $i$ with index lower than the two former pieces, the piece $u$ cannot be placed below piece $j$ if piece $j$ is placed on top of piece $i$. Thus, the constraints (\ref{ineq:new_inequality_identical_pieces}) below must be satisfied by all feasible solutions of the NFP-CM-VS model.
\begin{equation}
    \sum\limits_{k \in T_{ij}^{fg}} v_{ij}^{fgk} + \sum\limits_{k' \in B_{iu}^{fh}} v_{iu}^{fhk'} \leq 1, \qquad \forall (i,j,u,f,g,h) \in \{II_{iju}^{fgh} \:|\: t_j = t_u \}   \label{ineq:new_inequality_identical_pieces}
\end{equation}    

\subsubsection{Valid cuts and variable reductions between two pieces}

\begin{figure}[t]
\centering
\begin{minipage}[t]{0.45\textwidth}
\begin{tikzpicture}
\fill[black!5] (4.5,3) -- (4.5,4.5)  -- (3,4.5) -- (2.5,4) -- (2.5,3.5) -- (3,3) -- (4.5,3);
\draw[black] (4.5,3) -- (4.5,4.5)  -- (3,4.5) -- (2.5,4) -- (2.5,3.5) -- (3,3);
\draw[dashed, black] (3,3) -- (4.5,3);
\filldraw[black] (3,1.5) circle (0.05);
\draw (3.5,3.75) node {$A^1$};
\fill[black!5] (3,1.5) -- (4,1.5) -- (4.5,2) -- (4.5,3) -- (3,3) -- (2.5,2.5) -- (2.5,2) -- (3,1.5);
\draw[black] (3,3) -- (2.5,2.5) -- (2.5,2) -- (3,1.5);
\draw[black]  (4,1.5) -- (4.5,2) -- (4.5,3);
\draw[dashed, black] (3,1.5) -- (4,1.5);
\draw (3.5,2.25) node {$A^2$};
\fill[black!5] (3,0) -- (4,0) -- (4.5,0.5) -- (4.5,1) -- (4,1.5) -- (3,1.5) -- (2.5,1) -- (2.5,0.5) -- (3,0);
\draw[black] (3,0) -- (4,0) -- (4.5,0.5) -- (4.5,1) -- (4,1.5);
\draw[black] (3,1.5) -- (2.5,1) -- (2.5,0.5) -- (3,0);
\draw (3.5,0.75) node {$A^3$};
\fill[black!5] (5.5,2.5) -- (6.5,2.5) -- (6.5,3.5) -- (5.5,3.5) -- (5.5,2.5);
\draw[black] (5.5,2.5) -- (6.5,2.5) -- (6.5,3.5) -- (5.5,3.5) -- (5.5,2.5);
\filldraw[black] (5.5,2.5) circle (0.05);
\draw (6,3) node {$B^1$};
\draw (6,1.5) node {$A = A^1 \cup A^2 \cup A^3$};
\draw (5.75,1) node {$B = B^1$};
\draw (4.5,-0.5) node {(a) convex decomposition of pieces $A$ and $B$};
\end{tikzpicture}
\end{minipage}
\hfill
\begin{minipage}[t]{0.45\textwidth}
\begin{tikzpicture}[scale = 0.85]
\fill[black!5] (4.5,3) -- (4.5,4.5)  -- (3,4.5) -- (2.5,4) -- (2.5,3.5) -- (3,3) -- (4.5,3);
\draw (3.5,3.75) node {$A^1$};
\fill[black!5] (3,1.5) -- (4,1.5) -- (4.5,2) -- (4.5,3) -- (3,3) -- (2.5,2.5) -- (2.5,2) -- (3,1.5);
\draw (3.5,2.25) node {$A^2$};
\fill[black!5] (3,0) -- (4,0) -- (4.5,0.5) -- (4.5,1) -- (4,1.5) -- (3,1.5) -- (2.5,1) -- (2.5,0.5) -- (3,0);
\draw (3.5,0.75) node {$A^3$};
\draw[dashed, black] (2,2) -- (4.5,2)  -- (4.5,4.5) -- (2,4.5) -- (1.5,4) -- (1.5,2.5) -- (2,2);
\draw[dashed, black] (2,0.5) -- (4,0.5) -- (4.5,1) -- (4.5,3) -- (2,3) -- (1.5,2.5) -- (1.5,1) -- (2,0.5);
\draw[dashed, black] (2,-1) -- (4,-1) -- (4.5,-0.5) -- (4.5,1) -- (4,1.5) -- (2,1.5) -- (1.5,1) -- (1.5,-0.5) -- (2,-1);
\draw[black] (2.25,-0.25) node {$NFP_{AB}^{31}$};
\draw[black] (2.5,2.5) node {$NFP_{AB}^{21}$};
\draw[black] (2.25,3.75) node {$NFP_{AB}^{11}$};
\draw (3,-1.5) node {(b) convex NFP parts between pieces $A$ and $B$};
\end{tikzpicture}
\end{minipage}
\bigskip
\vspace{0.25cm}
\begin{minipage}[t]{0.3\textwidth}
\begin{tikzpicture}[scale = 0.85]
\fill[black!5] (1.5,-1) -- (1.5,4.5) -- (0,4.5) -- (0,-1) -- (1.5,-1);
\fill[black!5] (6,-1) -- (6,4.5) -- (4.5,4.5) -- (4.5,-1) -- (6,-1);
\draw[dashed, black] (2,2) -- (4.5,2)  -- (4.5,4.5) -- (2,4.5) -- (1.5,4) -- (1.5,2.5) -- (2,2);
\draw[dashed, black] (2,0.5) -- (4,0.5) -- (4.5,1) -- (4.5,3) -- (2,3) -- (1.5,2.5) -- (1.5,1) -- (2,0.5);
\draw[dashed, black] (2,-1) -- (4,-1) -- (4.5,-0.5) -- (4.5,1) -- (4,1.5) -- (2,1.5) -- (1.5,1) -- (1.5,-0.5) -- (2,-1);
\draw[black] (0.75,0.25) node {$v_{AB}^{31l}$};
\draw[black] (0.75,1.75) node {$v_{AB}^{21l}$};
\draw[black] (0.75,3.25) node {$v_{AB}^{11l}$};
\draw[black] (5.25,0.25) node {$v_{AB}^{31r}$};
\draw[black] (5.25,1.75) node {$v_{AB}^{21r}$};
\draw[black] (5.25,3.25) node {$v_{AB}^{11r}$};
\draw[black] (3,0) node {$NFP_{AB}^{31}$};
\draw[black] (3,2.5) node {$NFP_{AB}^{21}$};
\draw[black] (3,3.75) node {$NFP_{AB}^{11}$};
\draw (3,-1.5) node {(c) identical feasible left and};
\draw (3,-2) node   {right sub-regions in light grey};
\end{tikzpicture}
\end{minipage}
\hfill
\begin{minipage}[t]{0.2\textwidth}
\begin{tikzpicture}[scale = 0.85]
\fill[black!5] (2,2) -- (2,-2) -- (4.5,-2) -- (4.5,2) -- (2,2);
\fill[black!20] (2,-2) -- (4.5,-2) -- (4.5,-0.5) -- (4,-1) -- (2,-1) -- (2,-2);
\draw[black] (2,2) -- (4.5,2)  -- (4.5,4.5) -- (2,4.5) -- (1.5,4) -- (1.5,2.5) -- (2,2);
\draw[dashed, black] (2,-1) -- (4,-1) -- (4.5,-0.5) -- (4.5,1) -- (4,1.5) -- (2,1.5) -- (1.5,1) -- (1.5,-0.5) -- (2,-1);
\draw (3,0) node {$NFP_{AB}^{31}$};
\draw[black] (3,3.25) node {$NFP_{AB}^{11}$};
\draw[dashed, black]  (4,-2) -- (4,-1);
\draw[black] (3,1) node {$v_{AB}^{11b_1}$};
\draw[black] (3,-1.5) node {$v_{AB}^{31b_1}$};
\draw[black] (5.5,-1.5) node {$v_{AB}^{31b_2}$};
\draw[dashed, black] (4.25,-1.5) -- (5.1,-1.5);
\filldraw[black] (4.25,-1.5) circle (0.05);
\draw (3.5,-2.5) node {(d) subsumption between};
\draw (3.5,-3) node {regions of two convex parts};
\end{tikzpicture}
\end{minipage}
\hfill
\begin{minipage}[t]{0.3\textwidth}
\begin{tikzpicture}[scale = 0.85]
\fill[black!5] (1.5,2.5) -- (1.5,-2) -- (2,-2) -- (2,2) -- (1.5,2.5);
\fill[black!5] (4,-2) -- (4.5,-2) -- (4.5,-0.5) -- (4,-1) -- (4,-2);
\draw[black] (2,2) -- (4.5,2)  -- (4.5,4.5) -- (2,4.5) -- (1.5,4) -- (1.5,2.5) -- (2,2);
\draw[dashed, black] (2,-1) -- (4,-1) -- (4.5,-0.5) -- (4.5,1) -- (4,1.5) -- (2,1.5) -- (1.5,1) -- (1.5,-0.5) -- (2,-1);
\draw (3,0) node {$NFP_{AB}^{31}$};
\draw[black] (3,3.25) node {$NFP_{AB}^{11}$};
\draw[black] (0.5,1.75) node {$v_{AB}^{11b_0}$};
\draw[black] (5.5,-1.5) node {$v_{AB}^{31b_2}$};
\draw[dashed, black] (4.25,-1.5) -- (5.1,-1.5);
\filldraw[black] (4.25,-1.5) circle (0.05);
\draw[dashed, black] (0.95,1.75) -- (1.75,1.75);
\filldraw[black] (1.75,1.75) circle (0.05);
\draw (3.5,-2.5) node {(e) non-overlapping feasible};
\draw (3.5,-3) node {regions of two convex parts};
\end{tikzpicture}
\end{minipage}
\caption{This figure shows three geometric configurations between convex NFP parts of the pieces $A$ and $B$, which induces three types of logic relationships encoded as valid cuts and variable eliminations into our NFP-CM-VS models.}
\label{fig:cuts_between_two_convex_parts}
\end{figure}
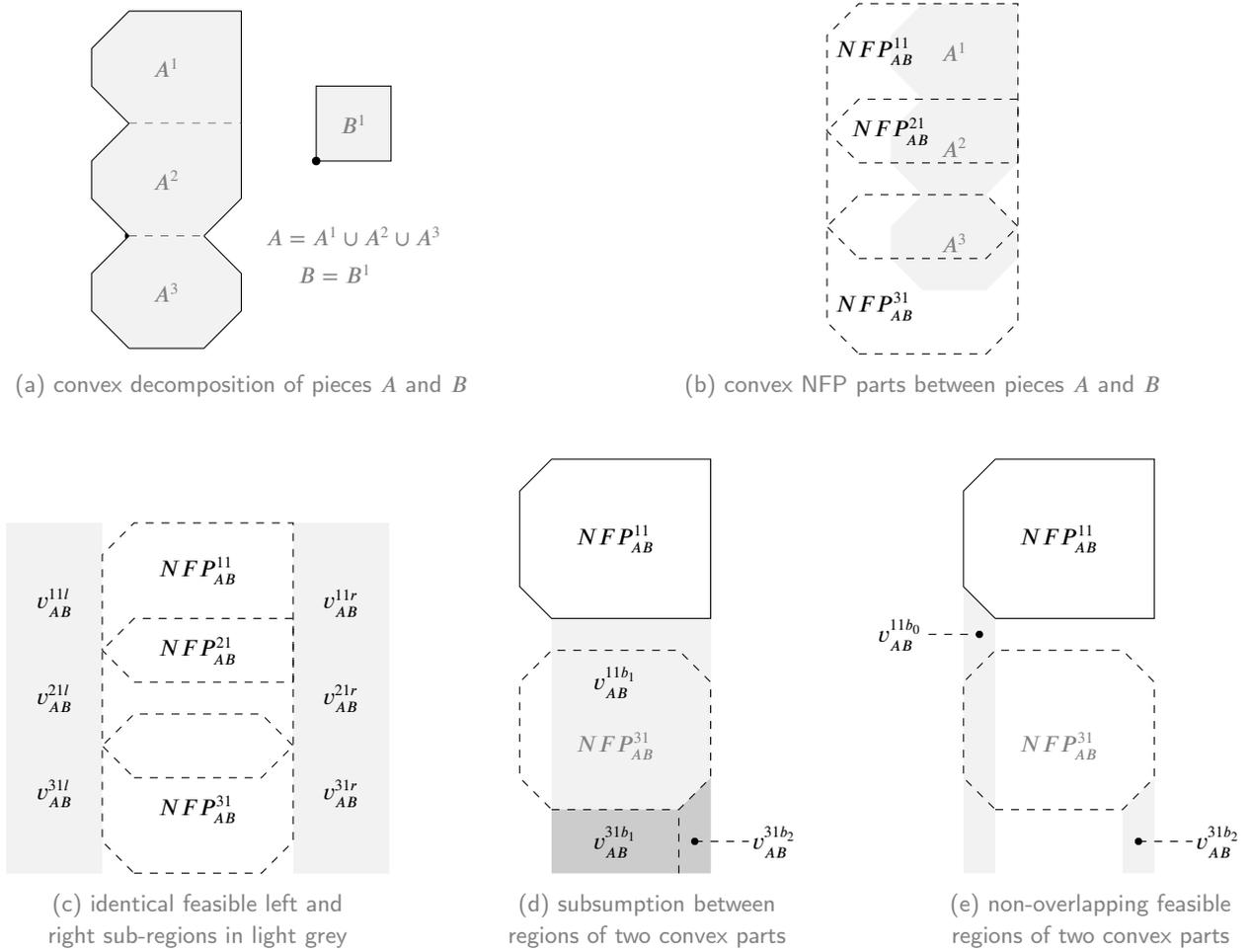

Given two non-convex pieces $i$ and $j$, their overall $NFP_{ij}$ is decomposed into a collection of convex $NFP_{ij}^{fg}$ parts obtained from the combination of the convex parts from both pieces, as shown in figure \ref{fig:cuts_between_two_convex_parts}.b. For example, figure \ref{fig:cuts_between_two_convex_parts}.a shows a convex piece $B$ and a non-convex piece $A$ decomposed into three convex parts, whilst figure \ref{fig:cuts_between_two_convex_parts}.b shows the three pairwise $NFP_{AB}^{fg}$ parts resulting from the combination of the convex parts of both former pieces, denoted by $NFP_{AB}^{11}$, $NFP_{AB}^{21}$, and $NFP_{AB}^{31}$. Each convex $NFP_{AB}^{fg}$ part in the example shown in figure \ref{fig:cuts_between_two_convex_parts} sets a collection of feasible sub-regions for the relative placement of piece $B$ regarding piece $A$. However, any feasible relative placement of the orbiting piece $B$ regarding piece $A$ in our MIP models must be in a common feasible region for all $NFP_{AB}^{fg}$ parts between both pieces. Thus, given two or more feasible sub-regions belonging to two or more different convex $NFP_{AB}^{fg}$ parts of $NFP_{AB}$, we identify three a priori geometric relationships between them inducing the set of logic relationships enumerated below, which allow tightening the formulation of our family of NFP-CM-VS models.

\begin{enumerate}
\item[(1)]\emph{Identical feasible sub-regions.} If feasible sub-regions of different convex NFP parts between two pieces represent the same feasible region for their relative placement, then we can eliminate all redundant binary variables and represent all identical regions by a single binary variable in the model. According to our convex decomposition in figure \ref{fig:convex_decomposition}, only the left and right feasible regions encoded by the binary variables $v_{ij}^{fgl}$ and $v_{ij}^{fgr}$ in different convex $NFP_{ij}^{fg}$ parts can be identical. For example, figure \ref{fig:cuts_between_two_convex_parts}.c shows that $v_{AB}^{11l}$, $v_{AB}^{21l}$, and $v_{AB}^{31l}$ encode the same feasible region on the left, whilst $v_{AB}^{11r}$, $v_{AB}^{21r}$, and $v_{AB}^{31r}$ encode the same feasible region on the right. Thus, we can set a single variable to enable the left and right feasible sub-regions and remove the remaining ones during the building of the model. This variable eliminations are defined by the equalities (\ref{eq:reduction1}) and (\ref{eq:reduction2}) below.
\begin{align}
v_{ij}^{fgl} &= v_{ij}^{f'g'l}, \qquad \forall (i,j,f,g,f',g') : (f,g) \neq (f',g') \land \underline{x}_{ij}^{f'g'} = \underline{x}_{ij}^{fg} \label{eq:reduction1} \\
v_{ij}^{fgr} &= v_{ij}^{f'g'r}, \qquad \forall (i,j,f,g,f',g') : (f,g) \neq (f',g') \land \overline{x}_{ij}^{f'g'} = \overline{x}_{ij}^{fg} \label{eq:reduction2}  
\end{align}
    
\item[(2)]\emph{Subsumed feasible sub-regions.} Figure \ref{fig:cuts_between_two_convex_parts}.d shows an example with two feasible sub-regions $R_{AB}^{31b_1}$ and $R_{AB}^{31b_2}$ respectively, which are subsumed by the feasible sub-region $R_{AB}^{11b_1}$ belonging to other convex NFP part between the same pieces, such that $R_{AB}^{31b_1} \subset R_{AB}^{11b_1}$ and $R_{AB}^{31b_2} \subset R_{AB}^{11b_1}$. Thus, this later subsumption relationships between feasible sub-regions set a logic implication between their corresponding binary variables, such that any feasible solution of the model that sets to 1 a binary variable enabling a subsumed sub-region must also set to 1 the binary variables enabling their subsumer regions from other convex NFP parts between the same pieces. Thus, we can include the valid inequalities (\ref{ineq:subsummed_regions}) in our models.
\begin{equation}
\sum_{k \in K_{ij}^{fg} \land (R_{ij}^{fgk} \subset K_{ij}^{f'g'k'})} v_{ij}^{fgk} \leq v_{ij}^{f'g'k'}, \qquad \forall (i,j,f,g,f',g') : (f,g) \neq (f',g')
\label{ineq:subsummed_regions}
\end{equation}
    
\item[(3)] \emph{Non-overlapping regions.} If two feasible sub-regions belonging to two different convex NFP parts between two pieces do not overlap, then their corresponding binary variables can be enabled in none feasible solution at the same time because the resulting MIP model would be infeasible. For example, figure \ref{fig:cuts_between_two_convex_parts}.e shows that $R_{AB}^{11b_0} \cap R_{AB}^{31b_2} = \emptyset$, which implies that any feasible solution of the model must satisfy that $v_{AB}^{11b_0} + v_{AB}^{31b_2} \leq 1$. Thus, the former logic relationships can be represented in the model by including the valid inequalities (\ref{ineq:non-overlapping_regions}).
\begin{equation}
v_{ij}^{fgk} + v_{ij}^{f'g'k'} \leq 1, \forall (i,j,f,g,f',g',k,k'): R_{ij}^{fgk} \cap R_{ij}^{f'g'k'} = \emptyset \label{ineq:non-overlapping_regions}
\end{equation}

\emph{Clique-based cuts.} The inequalities (\ref{ineq:non-overlapping_regions}) can be implemented in a much more efficient manner by defining them using the cliques of mutually-exclusive non-overlapping feasible sub-regions between convex NFP parts of two pieces, as defined by constraints (\ref{ineq:clique_based_cuts}). Given an undirected graph $G=(V,E)$, a set of vertexes $\mathcal{C} \subseteq V$ is a \emph{clique} if $\forall a,b \in \mathcal{C} \Rightarrow \exists (a,b) \in E$. The cliques allow factorizing many of the former constraints into a single one, which significantly reduces the number of constraints inserted into the model. Thus, we build an undirected graph $G = (V,E)$ for representing the collection of constraints (\ref{ineq:non-overlapping_regions}), such that $E$ contains an edge $(v_{ij}^{fgk}, v_{ij}^{f'g'k'})$ for each pair of binary variables satisfying (\ref{ineq:non-overlapping_regions}), and $V$ contains all the binary variables that take part in at least one of the former constraints. Next, we compute the cliques of non-overlapping feasible sub-regions using the minimum Edge Clique Covering (ECC) algorithm \cite{Conte2020-mo} provided by the ECC8 Java software library available at \url{https://github.com/Pronte/ECC}. Finally, the clique-based cuts (\ref{ineq:clique_based_cuts}) are inserted as SOS-type 1 constraints in the implementation of our models instead of constraints (\ref{ineq:non-overlapping_regions}).
\begin{equation}
\sum\limits_{v_{ij}^{fgk} \in \mathcal{C}_q} v_{ij}^{fgk} \leq 1, \qquad \forall \mathcal{C}_q \in \{\text{cliques covering constraints (\ref{ineq:non-overlapping_regions})}\} \label{ineq:clique_based_cuts}
\end{equation}
\end{enumerate}

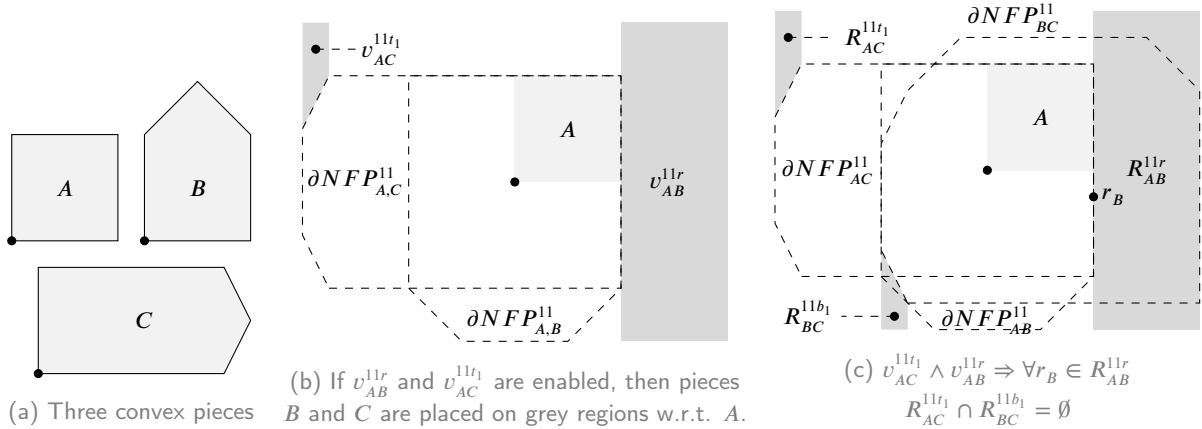
\begin{figure}[t!]
\begin{minipage}[t]{0.2\textwidth}
\begin{tikzpicture}[scale = 0.35]
\fill[black!5] (0,0) -- (4,0) -- (4,4) -- (0,4) -- (0,0);
\draw[black] (0,0) -- (4,0) -- (4,4) -- (0,4) -- (0,0);
\draw[black](2,2) node {$A$};
\filldraw[black] (0,0) circle (0.15);
\fill[black!5] (5,0) -- (9,0) -- (9,4) -- (7,6) -- (5,4) -- (5,0);
\draw[black] (5,0) -- (9,0) -- (9,4) -- (7,6) -- (5,4) -- (5,0);
\draw[black] (7,2) node {$B$};
\filldraw[black] (5,0) circle (0.15);
\fill[black!5] (1,-5) -- (8,-5) -- (9,-3) -- (8,-1) -- (1,-1) -- (1,-5);
\draw[black] (1,-5) -- (8,-5) -- (9,-3) -- (8,-1) -- (1,-1) -- (1,-5);
\draw[black] (5,-3) node {$C$};
\filldraw[black] (1,-5) circle (0.15);
\draw (4.5,-6.5) node {(a) Three convex pieces};
\end{tikzpicture}
\end{minipage}
\hfill
\begin{minipage}[t]{0.38\textwidth}
\begin{tikzpicture}[scale = 0.35]
\fill[black!5] (0,0) -- (4,0) -- (4,4) -- (0,4) -- (0,0);
\draw[black](2,2) node {$A$};
\filldraw[black] (0,0) circle (0.15);
\fill[black!15] (4,-6) -- (8,-6) -- (8,6) -- (4,6) -- (4,-6);
\draw[black](5.75,0) node {$v_{AB}^{11r}$};
\fill[black!15] (-7,4) -- (-7,6) -- (-8,6) -- (-8,2) -- (-7,4);
\draw[black](-5,5) node {$v_{AC}^{11t_1}$};
\draw[black, dashed] (-7.5,5) -- (-6,5);
\filldraw[black] (-7.5,5) circle (0.15);
\draw[black, dashed] (-4,-4) -- (-2,-6) -- (2,-6) -- (4,-4) -- (4,4) -- (-4,4) -- (-4,-4);
\draw[black, dashed] (-7,-4) -- (4,-4) -- (4,4) -- (-7,4) -- (-8,2) -- (-8,-2) -- (-7,-4);
\draw[black](0,-5.25) node {$\partial NFP_{A,B}^{11}$};
\draw[black](-6,0) node {$\partial NFP_{A,C}^{11}$};
\draw (0,-7.5) node {(b) If $v_{AB}^{11r}$ and $v_{AC}^{11t_1}$ are enabled, then pieces};
\draw (0,-8.75) node {$B$ and $C$ are placed on grey regions w.r.t. $A$.};
\end{tikzpicture}
\end{minipage}
\hfill
\begin{minipage}[t]{0.38\textwidth}
\begin{tikzpicture}[scale = 0.35]
\fill[black!5] (0,0) -- (4,0) -- (4,4) -- (0,4) -- (0,0);
\draw[black](2,2) node {$A$};
\filldraw[black] (0,0) circle (0.15);
\fill[black!15] (4,-6) -- (8,-6) -- (8,6) -- (4,6) -- (4,-6);
\draw[black](6,0) node {$R_{AB}^{11r}$};
\fill[black!15] (-7,4) -- (-7,6) -- (-8,6) -- (-8,2) -- (-7,4);
\draw[black](-4.5,5) node {$R_{AC}^{11t_1}$};
\draw[black, dashed] (-7.5,5) -- (-6,5);
\filldraw[black] (-7.5,5) circle (0.15);
\fill[black!15] (-3,-5) -- (-4,-3) -- (-4,-6) -- (-3,-6) -- (-3,-5);
\draw[black](-6.75,-5.5) node {$R_{BC}^{11b_1}$};
\draw[black, dashed] (-3.5,-5.5) -- (-5.5,-5.5);
\filldraw[black] (-3.5,-5.5) circle (0.15);
\draw[black, dashed] (-4,-4) -- (-2,-6) -- (2,-6) -- (4,-4) -- (4,4) -- (-4,4) -- (-4,-4);
\draw[black, dashed] (-7,-4) -- (4,-4) -- (4,4) -- (-7,4) -- (-8,2) -- (-8,-2) -- (-7,-4);
\draw[black](0,-5.65) node {$\partial NFP_{AB}^{11}$};
\draw[black](-6,0) node {$\partial NFP_{AC}^{11}$};
\draw[black, dashed] (-3,-5) -- (8,-5) -- (8,3) -- (6,5) -- (-1,5) -- (-3,3) -- (-4,1) -- (-4,-3) -- (-3,-5);
\filldraw[black] (4,-1) circle (0.15);
\draw[black](4.75,-1) node {$r_B$};
\draw[black](1,5.75) node {$\partial NFP_{BC}^{11}$};
\draw (0,-7.5) node {(c) $v_{AC}^{11t_1} \land v_{AB}^{11r} \Rightarrow \forall r_B \in R_{AB}^{11r}$};
\draw (0,-9) node {$R_{AC}^{11t_1} \cap R_{BC}^{11b_1} = \emptyset$};
\end{tikzpicture}
\end{minipage}
\caption{Figure (a) shows an example of three convex pieces whose relative placements produce a large number of unfeasible combinations of three binary variables. Figure (b) shows in grey the two feasible relative placements of pieces $B$ and $C$ regarding the static (pivot) piece $A$. Finally, figure (c) shows that if the binary variables $v_{AC}^{11t_1}$ and $v_{AB}^{11r}$ are enabled, then $\forall r_B \in R_{AB}^{11r}, R_{AC}^{11t_1} \cap R_{BC}^{11b_1} = \emptyset$. Thus, the activation of the relative placements enabled by $v_{BC}^{11b_1}$ will produce an unfeasible solution, as shown in figure (c). The joint activation of $v_{AC}^{11t_1}$ and $v_{AB}^{11r}$ above causes that the all relative placements of piece $C$ regarding piece $B$ be unfeasible, with the only exception of the left feasible sub-region enabled by $v_{BC}^{11r}$.}
\label{fig:unfeasibilty_cuts_triplets}
\end{figure}

\subsubsection{Valid inequalities among three pieces}

\emph{Unfeasible relative placements of three pieces}. One key advantage derived from our convex decomposition based on vertical slices is that it allows disclosing many unfeasible combinations of binary variables encoding feasible relative placements among multiple pieces, which can be explicitly inserted into the model as feasibility cuts. For instance, given three convex pieces $A$, $B$, and $C$ as shown in figure \ref{fig:unfeasibilty_cuts_triplets}.a, if binary variables $v_{AB}^{11r}$ and $v_{AC}^{11t_1}$ are enabled in any solution of the model, then pieces $B$ and $C$ will be placed in the grey sub-regions with respect to the piece $A$, as shown in figure \ref{fig:unfeasibilty_cuts_triplets}.b. However, figure \ref{fig:unfeasibilty_cuts_triplets}.c shows that if the three aforementioned pieces are placed in the former relative placements, then the piece $C$ can be placed in none relative feasible sub-region regarding piece $B$ that be compatible (feasible) with its relative placement regarding piece $A$ defined by $R_{AC}^{11t_1}$, with the only exception of the left feasible region enabled by $v_{BC}^{11l}$. In this way, many unfeasible relative placements among three pieces can be detected a priori and represented as valid inequalities in the model. Thus, we can identify many unfeasible relative placements among three pieces $(i,j,u)$ by testing if the relative placement of piece $u$ regarding piece $j$ is feasible in any solution of the problem given the relative placements of pieces $j$ and $u$ regarding piece $i$.

\emph{Feasibility cuts among multiple pieces}. Although we focus here on representing the unfeasible relative placements among three pieces, our ideas can be generalized to combinations of $(m \geq 3)$ pieces by fixing the pairwise relative placements of the first $m - 1$ pieces and testing the feasibility for the relative placements of the last piece regarding the first $m - 1$ pieces, at the expense of exponentially increasing the number of valid inequalities inserted into the model.

\begin{algorithm}[t!]
\caption{\emph{IsUnfeasible()} function uses interval arithmetic \cite[ch.2]{Moore2009-ut} to test if three fixed relative placements among three pieces produce an unfeasible solution by testing their potential overlapping along the horizontal direction. A closed interval of real numbers is denoted by $[a,b] = \{x \in \mathbb{R} :a \leq x \leq b\}$. Let be $X=[\underline{X},\overline{X}]$ and $Y=[\underline{Y},\overline{Y}]$ two intervals, we say that both are equal if they are the same subset of $\mathbb{R}$. The sum of two intervals is $Z = X + Y = [\underline{X} + \underline{Y},\overline{X} + \overline{Y}]$. The intersection of two intervals is $Z = X \cap Y = \{z: z \in X \land z \in Y\} = [\max\{\underline{X},\underline{Y}\},\min\{\overline{X},\overline{Y}\}]$, and  $X \cap Y = \emptyset$ if $\overline{Y} < \underline{X}$ or $\overline{X} < \underline{Y}$.}
\label{alg:three_pieces_test}
\begin{algorithmic}[1]
\Require $R_{ij}^{fgk}, R_{iu}^{fhk'}, R_{ju}^{ghk''} \subset \mathbb{R}^2$
\Ensure \emph{true} if $\nexists r_j \in R_{ij}^{fgk} | R_{iu}^{fhk'} \cap R_{ju}^{ghk''} \neq \emptyset$, or \emph{false} otherwise
\Function{IsUnfeasible}{$R_{ij}^{fgk}, R_{iu}^{fhk'}, R_{ju}^{ghk''}$}
\State $x_{ij}^{fgk} \gets [\min\{R_{ij,x}^{fgk}\}, \max\{R_{ij,x}^{fgk}\}] \subset \mathbb{R}$ \Comment{\text{sets the x-axis interval spanned by $R_{ij}^{fgk}$}}
\State $x_{iu}^{fhk'} \gets [\min\{R_{iu,x}^{fhk'}\}, \max\{R_{iu,x}^{fhk'}\}] \subset \mathbb{R}$ \Comment{\text{sets the x-axis interval spanned by $R_{iu}^{fhk'}$}}
\State $x_{ju}^{ghk''} \gets [\min\{R_{ju,x}^{ghk''}\}, \max\{R_{ju,x}^{ghk''}\}] \subset \mathbb{R}$ \Comment{\text{sets the x-axis interval spanned by $R_{ju}^{ghk''}$}}
\State $\delta x_{ju} \gets x_{ij}^{fgk} + x_{ju}^{ghk''} = [\underline{x}_{ij}^{fgk} + \underline{x}_{ju}^{ghk''}, \overline{x}_{ij}^{fgk} + \overline{x}_{ju}^{ghk''}]$ \Comment{sum of x-axis intervals}
\State $z \gets \delta x_{ju} \cap x_{iu}^{fhk'} = [\max\{\underline{\delta x}_{ju},\underline{x}_{iu}^{fhk'}\},\min\{\overline{\delta x}_{ju},\overline{x}_{iu}^{fhk'}\}]$   \Comment{gets feasible relative placement interval for $u$}
\State $r \in \{\text{true, false}\}$
\If{$z = \emptyset$} \Comment{feasibility test}
\State $\text{r} \gets \text{true}$ 
\Else
\State $\text{r} \gets \text{false}$ 
\EndIf
\State \Return $r$
\EndFunction
\end{algorithmic}
\end{algorithm}

\emph{Computation of feasibility cuts among three pieces}. Algorithm \ref{alg:three_pieces_test} partially evaluates the feasibility of $v_{ju}^{ghk''}$ given that $v_{ij}^{fgk}$ and $v_{iu}^{fhk'}$ are enabled by testing if $R_{ju}^{ghk''}$ and $R_{iu}^{fhk'}$ regions overlap along the horizontal direction using interval arithmetic \cite{Moore1966-cv, Moore2009-ut}. Then, all unfeasible relative placements among three pieces can be inserted into our model using combinatorial Benders cuts \cite{Codato2006-se}, as defined by constraints (\ref{ineq:unfeasibility_cuts_for_triplets}). Although our feasibility test only considers the overlapping along the horizontal direction, it is capable of detecting a large number of unfeasible relative placements among three pieces. Because the number of valid inequalities (\ref{ineq:unfeasibility_cuts_for_triplets}) grows rapidly, raising to tenths of millions for large problem instances, we insert this family of feasibility cuts into our NFP-CM-VS models as user cuts. Thus, the former valid inequalities are used by the Integer Programming (IP) solver to cut any feasible solution violating them, being inserted into the Branch and Cut (B\&C) exploration according to the heuristics rules implemented by each IP solver.
\begin{equation}
v_{ij}^{fgk} + v_{iu}^{fhk'} + \sum\limits_{k'' \in K_{ju}^{gh}} v_{ju}^{ghk''} \leq 2,  \forall (i,j,u,f,g,h,k,k',k''): v_{ij}^{fgk} + v_{iu}^{fhk'} = 2 \Rightarrow \forall r_B \in R_{ij}^{fgk}, R_{iu}^{fhk'} \cap R_{ju}^{ghk''} = \emptyset \label{ineq:unfeasibility_cuts_for_triplets}
\end{equation}

\subsection{The NFP-CM-VS2 model}
\label{sec:NFP_CM_VS2}

As mentioned above, our basic NFP-CM-VS model triples in average the number of constraints of the NFP-CM model \cite{Cherri2016-jf} to break the symmetries of the feasible space, unlike the Improved NFP-CM model \cite{Rodrigues2017-dy} that only doubles them. To solve this drawback, we introduce a reformulation of our NFP-CM-VS model, called NFP-CM-VS2, which is derived from the full NFP-CM-VS model defined by the objective function (\ref{model:nfp_cm_vs_objective}) and the constraints and variable eliminations (\ref{model:x_bounds}-\ref{ineq:unfeasibility_cuts_for_triplets}) by substituting the constraints (\ref{model:nfp_cm_vs_first_constraint}-\ref{model:nfp_cm_vs_last_constraint}) by the constraints (\ref{ineq:vs2_1}-\ref{ineq:vs2_2}) below.

\begin{align}
 & x_j - x_i \geq -(L_{ub} - (l_j^{min} + l_i^{max)}))v_{ij}^{fgl} + \overline{x}_{ij}^{fg}v_{ij}^{fgr} + \sum\limits_{k \in T_{ij}^{fg}} b_{ij,x}^{fgk}v_{ij}^{fgk} + \sum\limits_{k \in B_{ij}^{fg}} a_{ij,x}^{fgk}v_{ij}^{fgk}, \: \forall (i,j,f,g) \in I_{ij}^{fg} \label{ineq:vs2_1} \\
 & x_j - x_i \leq \underline{x}_{ij}^{fg}v_{ij}^{fgl} + (L_{ub} - (l_i^{min} + l_j^{max}))v_{ij}^{fgr} + \sum\limits_{k \in T_{ij}^{fg}} a_{ij,x}^{fgk}v_{ij}^{fgk} + \sum\limits_{k \in B_{ij}^{fg}} b_{ij,x}^{fgk}v_{ij}^{fgk}, \: \forall (i,j,f,g) \in I_{ij}^{fg} \label{ineq:vs2_2}
\end{align}

\emph{The most compact symmetry-breaking.} The constraints (\ref{ineq:vs2_1}) and (\ref{ineq:vs2_2}) introduce two significant improvements on the basic NFP-CM-VS model as follows. First, they remove the big-M factors from constraints (\ref{model:nfp_cm_vs_first_constraint}-\ref{model:nfp_cm_vs_last_constraint}). And second, they allow a significant reduction in the complexity of the model by factorizing all constraints encoding the vertical lines bounding the feasible sub-regions of each convex $NFP_{ij}^{fg}$ part into only two constraints, whilst they hold the symmetry-breaking of the symmetric solutions for the space of feasible relative placements between pieces with much fewer constraints per convex $NFP_{ij}^{fg}$ part than the state-of-the-art Improved NFP-CM model \cite{Rodrigues2017-dy}. Note that the Improved NFP-CM model requires 2 constraints per binary variable of each convex $NFP_{ij}^{fgk}$ part to cover the feasible space as shown in figure \ref{fig:convex_decomposition}.d, whilst the NFP-CM-VS2 model only requires 1 constraint per binary variable (ineq. \ref{ineq:main_edge_feasible_region}) plus 2 constraints (ineq. \ref{ineq:vs2_1}-\ref{ineq:vs2_2}) per convex $NFP_{ij}^{fgk}$ part. Thus, our NFP-CM-VS2 model provides the most compact symmetry-breaking of the feasible space of relative placements between pieces reported in the literature.

\section{Evaluation}
\label{sec:evaluatiomn}

The goals of the experiments in this section are as follows: (1) to evaluate the performance of our two new MIP models, called NFP-CM-VS and NFP-CM-VS2; (2) to evaluate the impact of the family of feasibility cuts among three pieces defined by constraints (\ref{ineq:unfeasibility_cuts_for_triplets}); (3) to carry out a fair comparison of the performance of our new models with the family of state-of-the-art continuous MIP models for irregular strip packing, which is based on the same hardware and software platform; (4) to replicate the state-of-the-art family of NFP-CM and Improved NFP-CM models from scratch; (5) to develop a reproducible benchmark of state-of-the-art MIP models based on our software implementation of all models evaluated herein into the same software library, which is provided as supplementary material (see Appendix \ref{sec:appendix_B}); (6) to evaluate the state-of-the-art Improved NFP-CM in a standard benchmark \cite{Cherri2016-jf} including many unexplored problem instances not considered in its introductory paper \cite{Rodrigues2017-dy}; (7) the independent confirmation of previous findings and results reported in the literature; and finally, (8) to elucidate the current state of the art on continuous MIP models for irregular strip packing in a sound and reproducible way.

\begin{table}[t!]
\centering
\begin{tabular}{lll}
MIP model & Reference & Implementation details \\
\hline
\\
NFP-CMnc  & \cite{Cherri2016-jf} & \makecell[l]{Exact replication of the NFP-CM model without cuts.} \\
\\
NFP-CM & \cite{Cherri2016-jf} & \makecell[l]{Exact replication of the NFP-CM model \cite{Cherri2016-jf} in which the non-overlapping feasibility \\ cuts \cite[Algo.1, step 14]{Cherri2016-jf} are inserted into the model as SOS-1 constraints for a fair \\ comparison with our models.} \\
\makecell[l]{Improved \\ NFP-CM \\ (baseline)} & \cite{Rodrigues2017-dy} & Exact replication of the Improved NFP-CM model as detailed in \cite[\S3-4]{Rodrigues2017-dy}.\\
NFP-CM-VSnc & this work & \makecell[l]{Objective function (\ref{model:nfp_cm_vs_objective}) with constraints (\ref{model:x_bounds}-\ref{ineq:clique_based_cuts}). All constraints are inserted into the \\ root node of the model, with the only exception of the clique-based cuts (\ref{ineq:clique_based_cuts}) that are \\ inserted into the model as SOS-1 constraints.} \\
\\
NFP-CM-VS & this work &  \makecell[l]{Same implementation than NFP-CM-VSnc plus the feasibility cuts among three \\ pieces (\ref{ineq:unfeasibility_cuts_for_triplets}) inserted as user cuts (Lazy = -1) into the model.}\\
\\
NFP-CM-VS2 & this work & \makecell[l]{Same implementation than NFP-CM-VS and substitution of the constraints (\ref{model:nfp_cm_vs_first_constraint}-\ref{model:nfp_cm_vs_last_constraint}) \\ by constraints (\ref{ineq:vs2_1}) and (\ref{ineq:vs2_2}), which are inserted into the root node of the model.}
\end{tabular}
\caption{Implementation details for all models evaluated herein. All models are built by ordering the pieces by non-increasing area.}
\label{tab:models_evaluated}
\end{table}

\subsection{Experimental setup}

We reproduce the same experiments carried-out by \citet{Cherri2016-jf} to evaluate the NFP-CMnc and NFP-CM models by replicating all models detailed in table \ref{tab:models_evaluated} and evaluate them in the same set of problem instances used by the former authors, with the exception of the instances with holes and rotations. Thus, our experiments include the thirty-five small problem instances and ten large ones shown in \citet[tables 2-3]{Cherri2016-jf}, which include most of problem instances evaluated by \citet{Alvarez-Valdes2013-wg} that were also considered by \citet{Cherri2016-jf}, plus one additional large instance. On the other hand, the state-of-the-art Improved NFP-CM model setting our baseline for comparison is evaluated here for the first time in most of small problem instances shown in table \ref{tab:results1} and all large problem instances shown in table \ref{tab:results2}.

Table \ref{tab:models_evaluated} details the state-of-the-art MIP models for irregular strip packing evaluated herein and the details of our software implementation. To study the impact of our new family of feasibility cuts among three pieces (\ref{ineq:unfeasibility_cuts_for_triplets}), we have defined two versions of our new NFP-CM-VS model to be evaluated in our experiments, which are called NFP-CM-VSnc and full NFP-CM-VS respectively, as shown in table \ref{tab:models_evaluated}. The only difference between these two later models is that NFP-CM-VSnc does not include our new family of feasibility cuts among three pieces (\ref{ineq:unfeasibility_cuts_for_triplets}).

All our experiments are based on our own software implementation of all MIP models evaluated herein into the same Java software library, called RAMNEST, which uses the Java API of Gurobi 9.5 to solve the models. Our source code and a pre-compiled version of our software are publicly available in our reproducibility dataset \cite{Lastra-Diaz2022-aw}. We have replicated the state-of-the-art family of NFP-CM models \cite{Cherri2016-jf, Rodrigues2017-dy} from scratch by integrating the three steps to evaluate all models as follows: (1) pre-processing to decompose the pieces into convex parts based on the CGAL implementation of the Greene's algorithm \cite{Greene1983-fl}; (2) in-memory building of the MIP models by calling the functions in the Java API for Gurobi 9.5; and (3) resolution of the optimization models by calling the Java API for Gurobi with its default parameters. As we mentioned above, we use the ECC8 Java software library \cite{Conte2020-mo} to compute the cliques used to define the clique-based cuts (\ref{ineq:clique_based_cuts}) inserted as SOS-1 constraints in our NFP-CM-VS and NFP-CM-VS2 models. The version of the ECC algorithm \cite{Conte2020-mo} provided by ECC8 is not deterministic, which means that subsequent calls to ECC8 with the same data could return different cliques. For this reason, we save the cliques computed by the ECC8 software library into plain text files to allow the exact replication of our experiments, and to study isolately the impact of the valid cuts (\ref{ineq:unfeasibility_cuts_for_triplets}) in our two former models regarding the NFP-CM-VSnc model. All our experiments were implemented on a desktop UBUNTU 20.04 computer with an AMD Ryzen 7 5800X @3.8 GHz CPU (8 cores) and 64 Gb RAM.

\subsection{Reproducing our benchmarks}

All our experiments were generated by running the \emph{RamNestingdriver} program distributed with RAMNEST V1R1 software library \cite{Lastra-Diaz2022-aw} with two reproducible benchmark files that defines all problem instances evaluated herein. All problem instances are defined in ESICUP XML file format and provided with our raw output data in our reproducibility dataset \cite{Lastra-Diaz2022-aw}. The evaluation of each MIP model generates a Scalable Vector Graphics (SVG) file providing an image of the solution, and a raw output file in (*.csv) file format reporting the following data for each problem instance (see  Appendix \ref{sec:appendix_A}): (1) name of the problem instance; (2) number of pieces; (3) nesting efficiency; (4-5) lower and upper bounds of the solution; (6) MIP gap; (7) number of binary variables; (8) number of B\&B nodes; (9) number of simplex iterations; and (10) running time in seconds. All data reported herein is automatically generated by running the \emph{benchmark\_results} R-language script file on the collection of raw output files. Finally, all our models, experiments, and results can be exactly reproduced by following the instructions detailed in Appendix B, which also explains how to evaluate the MIP models detailed in table \ref{tab:models_evaluated} in other problem instances provided in ESICUP XML file format.

\subsection{Evaluation metrics}

To evaluate the performance of the MIP models detailed in table \ref{tab:models_evaluated}, we compare the upper bound, solution gap $(\frac{UB - LB}{UB})$, and computation time in seconds obtained by each MIP model in the evaluation of all problem instances. In order to provide a fair and unbiased comparison of the performance of all MIP models evaluated herein, we adopt the same approach proposed by \citet{Cherri2016-jf}, which uses performance profiles \cite{Dolan2002-wd} based on the ratios between the computation time of each model and the best computation time obtained by any model as a performance metric. Let be $\Phi$ and $M$ the sets of problem instances and MIP models evaluated herein, respectively. Then the ratio of computation times $t_{\phi,m}$  for each model $m \in M$ is defined by $r_{\phi,m}$ in formula (\ref{eq:ratio}), considering that the computation time is infinite whenever the models are unable to solve the problem instance up to the optimality. We assume an arbitrary parameter $r_M \geq r_{\phi,m} \: \forall (\phi,m) \in \Phi \times M$, such that $r_{\phi,m} = r_M$ if and only if the model $m$ is unable to solve the problem $\phi$. \citet{Dolan2002-wd} define a \emph{performance profile} of a solver or optimization model as the Cumulative Distribution Function (CDF) of the computation time ratios $r_{\phi,m}$, denoted by $\rho_m(\tau)$, and defined in formula (\ref{eq:performance_profile}) below. The performance profiles defined by the $\rho_m(\tau)$ function set a well-founded and broadly adopted metric to compare optimization models by avoiding any bias derived from a particular set of problem instances and dealing with those cases in which the models are unable to solve the problem up to the optimality.
\begin{align}
r_{\phi,m} &= \frac{t_{\phi,m}}{\min \{t_{\phi,i}: i \in M \}}, \quad \forall (\phi, m) \in \Phi \times M \label{eq:ratio} \\
\rho_m(\tau) &= \frac{1}{n_\phi} \{\phi \in \Phi : r_{\phi,m} \leq \tau\}, \quad \forall (m,\tau) \in M \times [1, r_M] \label{eq:performance_profile}
\end{align}
\subsection{Results obtained}

\begin{figure}[t]
\centering
\includegraphics[scale=0.55]{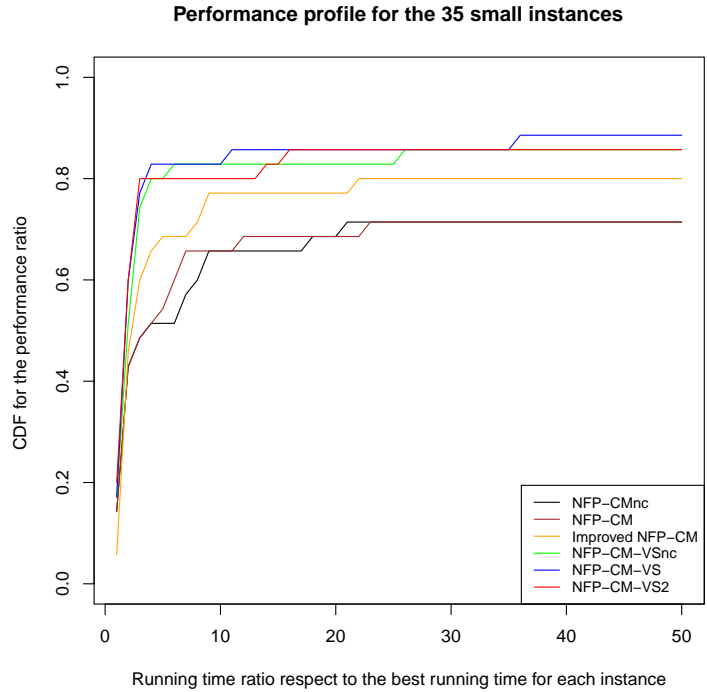}
\caption{Performance profile \cite{Dolan2002-wd} showing the Cumulative Distribution Function (CDF) for the performance ratio $r_{\phi,m}$ comparing the running times of all MIP models in the evaluation of the small problem instances.}
\label{fig:profile_small_instances}
\end{figure}

Table \ref{tab:results1} reports the terminating gap and running time in seconds obtained by the MIP models in the evaluation of the thirty-five small problem instances, whilst table \ref{tab:results2} reports the terminating gap for the eleven large instances. Figure \ref{fig:profile_small_instances} shows the performance profile curves comparing the performance ratio (\ref{eq:performance_profile}) obtained by all models in the evaluation of the small instances. Because of the lack of room, tables \ref{tab:results1} and \ref{tab:results2} omit the presentation of the lower and upper bounds. However, this missing information and all our raw output data are provided in Appendix A as supplementary material.

\begin{sidewaystable}[h]
\caption{Terminating gap and running time in seconds obtained by the MIP models in the evaluation of small problem instances with \# pieces. All models were implemented into the same software library based on Gurobi 9.5 and UBUNTU 20.04. Time Limit (TL) was set to 3600 seconds. Best results are shown in bold.}
\label{tab:results1}
\begin{tabular}{lccccccc|cccccc}
\multicolumn{2}{c}{} & \multicolumn{2}{c}{NFP-CMnc \cite{Cherri2016-jf}} &
\multicolumn{2}{c}{NFP-CM \cite{Cherri2016-jf}} &
\multicolumn{2}{c}{Improved NFP-CM \cite{Rodrigues2017-dy}} &
\multicolumn{2}{c}{\small{NFP-CM-VSnc}} &
\multicolumn{2}{c}{\small{NFP-CM-VS}} &
\multicolumn{2}{c}{\small{NFP-CM-VS2}} \\
\small{Instance} & \small{\#} & \small{GAP} & \small{\makecell{Time \\ (secs)}} & \small{GAP} & \small{\makecell{Time \\ (secs)}} & \small{GAP} & \small{\makecell{Time \\ (secs)}} & \small{GAP} & \small{\makecell{Time \\ (secs)}} & \small{GAP} & \small{\makecell{Time \\ (secs)}} & \small{GAP} & \small{\makecell{Time \\ (secs)}} \\
\hline
three & 3 & 0 & 0.005 & 0 & 0.006 & 0 & 0.008 & 0 & \textbf{ 0.004 } & 0 & 0.005 & 0 & 0.005 \\ 
  Shapes4 & 4 & 0 & 0.162 & 0 & 0.42 & 0 & \textbf{ 0.138 } & 0 & 0.306 & 0 & 0.272 & 0 & 0.304 \\ 
  glass1 & 5 & 0 & 0.037 & 0 & \textbf{ 0.016 } & 0 & 0.024 & 0 & 0.054 & 0 & 0.05 & 0 & 0.035 \\ 
  fu5 & 5 & 0 & 0.028 & 0 & 0.029 & 0 & 0.024 & 0 & 0.020 & 0 & \textbf{0.019} & 0 & 0.021 \\ 
  fu6 & 6 & 0 & \textbf{ 0.033 } & 0 & \textbf{ 0.033 } & 0 & 0.039 & 0 & 0.052 & 0 & 0.066 & 0 & 0.04 \\ 
  threep2 & 6 & 0 & \textbf{ 0.47 } & 0 & 0.474 & 0 & 0.54 & 0 & 1.255 & 0 & 0.971 & 0 & 1.242 \\ 
  threep2w9 & 6 & 0 & 0.658 & 0 & 0.668 & 0 & 0.647 & 0 & \textbf{ 0.528 } & 0 & 0.595 & 0 & 0.553 \\ 
  fu7 & 7 & 0 & \textbf{ 0.087 } & 0 & 0.09 & 0 & 0.115 & 0 & 0.145 & 0 & 0.092 & 0 & 0.146 \\ 
  glass2 & 7 & 0 & 0.086 & 0 & \textbf{ 0.064 } & 0 & 1.386 & 0 & 5.306 & 0 & 2.292 & 0 & 0.981 \\ 
  Shapes8 & 8 & 0.115 & TL & 0.115 & TL & 0.094 & TL & 0 & 288.602 & 0 & 111.896 & 0 & \textbf{ 109.222 } \\ 
  fu8 & 8 & 0 & 0.711 & 0 & \textbf{ 0.71 } & 0 & 0.751 & 0 & 0.812 & 0 & 0.749 & 0 & 0.963 \\ 
  fu9 & 9 & 0 & \textbf{ 0.784 } & 0 & 0.803 & 0 & 1.482 & 0 & 1.654 & 0 & 0.982 & 0 & 1.458 \\ 
  glass3 & 9 & 0 & \textbf{ 1.749 } & 0 & 2.366 & 0 & 412.119 & 0 & 44.624 & 0 & 18.862 & 0 & 23.21 \\ 
  threep3 & 9 & 0 & 191.169 & 0 & 190.446 & 0 & 274.271 & 0 & 183.403 & 0 & 206.237 & 0 & \textbf{ 114.837 } \\ 
  threep3w9 & 9 & 0.045 & TL & 0.045 & TL & 0 & 305.582 & 0 & 282.139 & 0 & 663.228 & 0 & \textbf{ 244.758 } \\ 
  Dighe2 & 10 & 0 & \textbf{ 3.227 } & 0 & 3.399 & 0 & 8.341 & 0 & 8.801 & 0 & 10.365 & 0 & 6.399 \\ 
  fu10 & 10 & 0 & 120.169 & 0 & 120.88 & 0 & \textbf{ 17.724 } & 0 & 34.216 & 0 & 21.944 & 0 & 31.438 \\ 
  J1-10-20-0 & 10 & 0 & 7.333 & 0 & 4.728 & 0 & 6.047 & 0 & 4.203 & 0 & 4.998 & 0 & \textbf{ 3.341 } \\ 
  J1-10-20-1 & 10 & 0 & 181.061 & 0 & 48.869 & 0 & 18.751 & 0 & 27.881 & 0 & 18.669 & 0 & \textbf{ 8.886 } \\ 
  J1-10-20-2 & 10 & 0 & 1.458 & 0 & \textbf{ 1.407 } & 0 & 2.997 & 0 & 3.475 & 0 & 4.005 & 0 & 3.037 \\ 
  J1-10-20-3 & 10 & 0 & 3591.792 & 0 & 729.372 & 0 & 110.605 & 0 & 35.642 & 0 & \textbf{ 33.049 } & 0 & 88.985 \\ 
  J1-10-20-4 & 10 & 0 & 1362.745 & 0 & 1884.525 & 0 & 26.273 & 0 & 24.6 & 0 & 29 & 0 & \textbf{ 20.173 } \\ 
  fu & 12 & 0 & 1047.891 & 0 & 1047 & 0 & 248.727 & 0 & \textbf{ 149.868 } & 0 & 256.79 & 0 & 336.361 \\ 
  J1-12-20-0 & 12 & 0 & 2683.008 & 0.083 & TL & 0 & 169.103 & 0 & 91.278 & 0 & \textbf{ 36.777 } & 0 & 43.46 \\ 
  J1-12-20-1 & 12 & 0 & 92.233 & 0 & 133.258 & 0 & 70.625 & 0 & 58.995 & 0 & 24.367 & 0 & \textbf{ 23.723 } \\ 
  J1-12-20-2 & 12 & 0 & 211.185 & 0 & 74.309 & 0 & 207.579 & 0 & 33.762 & 0 & \textbf{ 27.133 } & 0 & 61.211 \\ 
  J1-12-20-3 & 12 & 0 & 1774.853 & 0 & 944.026 & 0 & 381.156 & 0 & \textbf{ 215.248 } & 0 & 302.673 & 0 & 259.213 \\ 
  J1-12-20-4 & 12 & 0.154 & TL & 0.077 & TL & 0 & 1290.522 & 0 & 161.757 & 0 & \textbf{ 151.04 } & 0 & 189.593 \\ 
  J1-14-20-0 & 14 & 0 & 3404.836 & 0.083 & TL & 0 & 1554.604 & 0 & \textbf{ 191.53 } & 0 & 394.599 & 0 & 257.842 \\ 
  J1-14-20-1 & 14 & 0.167 & TL & 0.118 & TL & \textbf{ 0.029 } & TL & \textbf{ 0.029 } & TL & 0.118 & TL & \textbf{ 0.029 } & TL \\ 
  J1-14-20-2 & 14 & \textbf{ 0.143 } & TL & \textbf{ 0.143 } & TL & \textbf{ 0.143 } & TL & \textbf{ 0.143 } & TL & \textbf{ 0.143 } & TL & \textbf{ 0.143 } & TL \\ 
  J1-14-20-3 & 14 & 0 & 786.297 & 0 & 1067.409 & 0 & 345.419 & 0 & \textbf{ 92.546 } & 0 & 221.852 & 0 & 115.572 \\ 
  J1-14-20-4 & 14 & \textbf{ 0.143 } & TL & \textbf{ 0.143 } & TL & \textbf{ 0.143 } & TL & \textbf{ 0.143 } & TL & \textbf{ 0.143 } & TL & \textbf{ 0.143 } & TL \\ 
  Poly1a & 15 & \textbf{ 0.153 } & TL & 0.156 & TL & 0.196 & TL & 0.172 & TL & 0.212 & TL & 0.173 & TL \\ 
  Dighe1 & 16 & 0 & 587.735 & 0 & 1009.445 & 0.224 & TL & 0 & 2208.002 & 0 & \textbf{ 380.744 } & 0.197 & TL
\end{tabular}
\end{sidewaystable}
\clearpage

\begin{table}[h]
\caption{Terminating gap $(\frac{UB - LB}{UB})$ and running time in seconds obtained by the MIP models in the evaluation of all large problem instances. All models were implemented into the same software library based on Gurobi 9.5 onto an UBUNTU 20.04 computer. Time Limit (TL) was set to 3600 seconds, but none model was able to solve any instance up to optimality within this time limit. The 'X' symbol denotes that no feasible solution was obtained. Best gap results are shown in bold.}
\label{tab:results2}
\begin{tabular}{lccccccc}
 & & \small{NFP-CMnc} & \small{NFP-CM} & \small{\makecell{Improved NFP-CM}} &
\small{NFP-CM-VSnc} & \small{NFP-CM-VS} & \small{NFP-CM-VS2} \\
\small{Instance} & \small{\#} & \small{GAP} & \small{GAP} & \small{GAP} & \small{GAP} & \small{GAP} & \small{GAP} \\
\hline
Blaz2-16 & 16 & 0.299 & \textbf{ 0.285 } & 0.322 & 0.304 & 0.302 & 0.312 \\ 
Blaz2 & 20 & \textbf{ 0.265 } & 0.271 & 0.296 & 0.269 & \textbf{ 0.265 } & 0.281 \\ 
Mao & 20 & \textbf{ 0.243 } & 0.285 & X & 0.798 & 0.297 & X \\ 
Albano & 24 & 0.235 & 0.241 & X & 0.259 & \textbf{ 0.19 } & 0.25 \\ 
Marques & 24 & \textbf{ 0.176 } & 0.207 & X & X & 0.188 & 0.179 \\ 
Jakobs1 & 25 & \textbf{ 0.167 } & \textbf{ 0.167 } & \textbf{ 0.167 } & \textbf{ 0.167 } & \textbf{ 0.167 } & \textbf{ 0.167 } \\ 
Jakobs2 & 25 & 0.321 & 0.296 & 0.367 & 0.321 & 0.321 & \textbf{ 0.269 } \\ 
Blaz1 & 28 & 0.243 & 0.236 & X & 0.262 & \textbf{ 0.232 } & 0.264 \\ 
Dagli & 30 & 0.278 & 0.249 & X & X & 0.27 & \textbf{ 0.215 } \\ 
Shapes0 & 48 & X &  X & X & X & X & X \\ 
Trousers & 64 & \textbf{ 0.235 } & 0.241 & X & X & X & X 
\end{tabular}
\end{table}

\section{Discussion}
\label{sec_discussion}

\subsection{Comparison in small problem instances}

The entire family of NFP-CM-VS models obtains a significantly higher performance ratio than the family of NFP-CM models, and significantly outperforms the baseline Improved NFP-CM model, as shown in figure \ref{fig:profile_small_instances}. Likewise, the NFP-CM-VSnc and NFP-CM-VS models solve the largest number of small problem instances up to optimality, as shown in table \ref{tab:results1}.

The symmetry-breaking proposed by our new geometric decomposition based on vertical slices significantly outperforms that proposed by the baseline Improved NFP-CM model, as shown in figure \ref{fig:profile_small_instances}, and the results reported for large problem instances in table \ref{tab:results2}, in which the NFP-CM-VS and NFP-CM-VS2 models obtain better terminating gaps than the Improved NFP-CM model in all large instances, moreover doubling the number of instances in which a feasible solution is obtained.

Looking at the results in table \ref{tab:results1}, we find that among the small instances solved up to optimality, NFP-CM-VS2 is the fastest model in 7 instances, whilst NFP-CM-VS, NFP-VSnc, and NFP-CMnc models are the fastest ones in 6 instances, NFP-CM in 5, and Improved NFP-CM in 2.

The NFP-CM-VS and NFP-CM-VSnc models obtain a slightly better performance ratio than the NFP-CM-VS2 model, as shown in figure \ref{fig:profile_small_instances}. Thus, the compact formulation of the NFP-CM-VS2 model is unable to outperform the NFP-CM-VS models. Despite our NFP-CM-VS2 model defining a much more compact formulation than the NFP-CM-VS models and removing almost two-thirds of constraints and big-M terms, not only does NFP-CM-VS2 not improve on the NFP-CM-VS models, but it also performs slightly worse than the later ones. Thus, we conjecture that the NFP-CM-VS formulation provides a more tightened model than the more compact NFP-CM-VS2 formulation, making the B\&C exploration of the MIP solver faster.

\subsection{Comparison in large problem instances}

The entire family of NFP-CM-VS models obtains lower terminating gaps than the baseline Improved NFP-CM model in the evaluation of large problem instances. Moreover, our new family of MIP models doubles the cases in which it can find feasible solutions, as shown in table \ref{tab:results2}.

Current state-of-the-art MIP models evaluated herein are unable to either solve a large instance up to optimality or obtain feasible solutions for all of them, as shown in table \ref{tab:results2}. NFP-CMnc and NFP-CM models obtain feasible solutions for 10 of 11 large problem instances, whilst NFP-CM-VS does it for 9 cases, NFP-CM-V2 for 8, NFP-CM-VSnc for 7, and the Improved NFP-CM for only 4.

The NFP-CMnc obtains the lowest terminating gap in 4 large problem instances, whilst NFP-CM-VS does it in 3 cases, NFP-CM-VS2 in 2, and NFP-CM in 1, as shown in table \ref{tab:results2}.

NFP-CM and NFP-CMnc models significantly outperform the baseline Improved NFP-CM model in the large problem instances, as shown in table \ref{tab:results2}. Table \ref{tab:results2} also shows that Improved NFP-CM fails in obtaining feasible solutions in most of instances and it outperform NFP-CM model in none large instance. Thus, the symmetry-breaking proposed by the Improved NFP-CM is unable to provide consistent performance results in the large problems instances in comparison with its results in the small ones.

NFP-CMnc obtains competitive performance results in the large problem instances regarding NFP-CM, despite it have not been evaluated before, as shown in table \ref{tab:results2}. This conclusion is relevant because despite the introductory paper of \citet{Cherri2016-jf} omitted the evaluation of NFP-CMnc in the large problem instances, our results show that there is no a significant difference in performance between NFP-CMnc and NFP-CM, both in small and large problem instances, as shown in figure \ref{fig:profile_small_instances} and table \ref{tab:results2}.

\subsection{Impact of valid cuts among three pieces}

The valid cuts among three pieces implemented by our NFP-CM-VS model significantly improve the performance of the NFP-CM-VSnc model. This conclusion can be drawn by comparing the performance profiles shown in figure \ref{fig:profile_small_instances} and the termination gap values reported in table \ref{tab:results2} for both former models. Looking at the columns of both aforementioned models in table \ref{tab:results2}, you can see that NFP-CM-VS obtains a lower or equal terminating gap than NFP-CM-VSnc in all instances in which at least a feasible solution is found. Moreover, NFP-CM-VS finds feasible solutions for two large instances more than NFP-CM-VSnc. Thus, these performance improvements can only be attributed to the single difference between both former models: the valid cuts among three pieces defined by inequalities (\ref{ineq:unfeasibility_cuts_for_triplets}).

\subsection{Impact of technological advances and confirmation of previous results}

Our replication of all MIP models evaluated herein allows confirming previous results and drawing some conclusions on the impact of the technological advances both in hardware and MIP solvers regarding the results reported by \citet{Cherri2016-jf} and \citet{Rodrigues2017-dy}, regardless our experiments are based on Gurobi 9.5 (2022) instead of the CPLEX 12.6 solver used by the former authors.

The Improved NFP-CM baseline model obtains a higher performance ratio than the NFP-CM and NFP-CMnc models in the evaluation of small problem instances, as shown in figure \ref{fig:profile_small_instances}. First, this conclusion allows confirming in a sound way that the Improved NFP-CM sets the current state-of-the-art of the problem, and thus, it sets our baseline for comparison. We note that the Improved NFP-CM model introduced by \citet{Rodrigues2017-dy} was not evaluated in the same set of problem instances than the NFP-CM model \cite{Cherri2016-jf}, such as done here for the first time. On the contrary, the Improved NFP-CM model was evaluated in a different and more reduced set of instances, which limits the scope of previous conclusions. However, we bridge this minor gap here by independently confirming the achievements of the family of NFP-CM models \cite{Cherri2016-jf, Rodrigues2017-dy} in a more sound way.

The Dighe1 instance is solved for the first time up to the optimality by the NFP-CMnc, NFP-CM, NFP-CM-VSnc, and NFP-CM-VS models. This conclusion can be drawn by looking the results reported for the \emph{Dighe1} instance in table \ref{tab:results1} and comparing them with the results reported in \cite[table 1]{Cherri2016-jf}.

The advances in hardware and MIP solvers contribute to improve the performance of NFP-CMnc and NFP-CM models. This conclusion can be drawn by looking the results reported in tables \ref{tab:results1} and \ref{tab:results2}, and comparing them with those reported in \cite[tables 1-2]{Cherri2016-jf}. For instance, NFP-CMnc and NFP-CM are able to solve for the first time the \emph{threep3}, \emph{J1-10-20-3}, \emph{fu}, and \emph{Dighe1} instances. However, NFP-CM cannot solve the \emph{Shapes8} and \emph{threep3w9} instances up to the optimality in our experiments, which we attribute to differences in B\&C algorithms of the MIP solvers and the implementation of some valid inequalities using SOS-1 type constraints, as detailed in table \ref{tab:models_evaluated}. On the other hand, this technology improvement is much more noticeable in the case of the large problem instances reported in table \ref{tab:results2}, in which NFP-CM significantly reduce the terminating gap in all instances and it is able to obtain for the first time a feasible solution for the \emph{Trousers} instance.

\subsection{The new state of the art}

Our NFP-CM-VS model sets the new state-of-the-art among the family of continuous MIP models for nesting and provides a more consistent performance ratio than current state-of-the-art models, as shown in figure \ref{fig:profile_small_instances} and tables \ref{tab:results1} and \ref{tab:results2}. Although NFP-CM-VS only obtains comparable results to that obtained by the NFP-CMnc and NFP-CM models in terms of terminating gap in the evaluation of large problem instances, in which NFP-CMnc has the advantage of having one third less constraints than NFP-CM-VS, which significantly reduces the resolution of LP problems during the B\&C exploration, the performance ratios reported in figure \ref{fig:profile_small_instances} show that all improvements proposed in the formulation of our NFP-CM-VS models improve the state-of-the-art and suggest new lines of improvement for the problem. Moreover, NFP-CM-VS shows consistent performance results in the full range of problem instances evaluated herein, unlike the Improved NFP-CM model whose performance significantly decreases in the large problem instances.

\section{Conclusions and future work}
\label{sec:conclusions}

We have introduced a new family of continuous MIP models for irregular strip packing with two different formulations, abbreviated NFP-CM-VS and NFP-CM-VS2, which is based on a new convex decomposition of the feasible regions between convex parts into vertical slices, together with a new family of valid inequalities, symmetry breakings, and variable eliminations derived from the former geometric decomposition.

Our new family of MIP models outperform the state-of-the-art family of NFP-CM models introduced by \citet{Cherri2016-jf} and \citet{Rodrigues2017-dy}. Our new NFP-CM-VS model significantly and consistently sets the new state of the art of the problem. We show that our new geometric decomposition based on vertical slices outperforms the symmetry-breaking proposed by the Improved NFP-CM model \cite{Rodrigues2017-dy}, despite NFP-CM-VS tripling the number of constraints of NFP-CM \cite{Cherri2016-jf} instead of only doubling it as done by the Improved NFP-CM model. Likewise, we show that our new family of valid cuts among three pieces significantly contribute to the performance gain of our NFP-CM-VS model.

Another significant contribution is the introduction of the first confirmatory and reproducible experimental survey in this line of research, which is based on our software implementation of all models evaluated herein into a same Java software library, together with a detailed reproducibility protocol and dataset provided as supplementary material to allow the exact replication of all our models, experiments, and results.

We confirm the previous achievements of the family of NFP-CM models \cite{Cherri2016-jf}, setting the outperformance of the Improved NFP-CM model \cite{Rodrigues2017-dy} on the former ones in a sound way by replicating all models from scratch and comparing them in the same benchmark reported by \citet{Cherri2016-jf}.

As forthcoming activities, we plan to study the proposal of any decomposition method or matheuristics that allows improve the performance of our new family of MIP models for nesting.

\section*{Acknowledgements}

We are grateful to Luiz Cherri, Ramón Álvarez-Valdés, and Antonio Martínez-Sykora for providing us the problem instances used in their benchmarks, and to Alicia Lara-Clares by testing our reproducibility protocol.

\appendix

\section{Appendix A: Raw output data for all models}
\label{sec:appendix_A} 

This appendix introduces all raw output data generated by our experiments not reported herein by lack of room.

\section{Appendix B: The reproducible experiments on irregular strip-packing}
\label{sec:appendix_B}

This appendix introduces a detailed reproducibility protocol and dataset \cite{Lastra-Diaz2022-aw} providing our raw output data together with a collection of software and data resources to allow the exact replication of all our experiments and results.
\printcredits

\bibliographystyle{cas-model2-names}

\bibliography{biblio.bib}


\end{document}